\documentclass[10pt]{amsart}
\usepackage{amssymb,amsmath,amsfonts,amssymb,mathtools}
\usepackage{latexsym,euscript,exscale}
\usepackage[english]{babel}
\usepackage{graphicx}
\usepackage{hyperref}
\usepackage{float}
\usepackage{xcolor}
\newcommand{\pf}{
\noindent {\it Proof: }}

\newcommand{\D}{\mathbb{D}}

\newcommand{\Finf}{{\mathcal F}^{\infty}}

\newcommand {\N}{\mathbb N}

\newcommand {\R}{\mathbb R}
\newcommand {\Z}{\mathbb Z}
\newcommand {\C}{\mathbb C}

\newcommand{\pff}{$\hfill \square$}
\newcommand{\norm}[1]{\ensuremath{\left\|#1\right\|}}
\newcommand{\abs}[1]{\ensuremath{\left|#1\right|}}
\newtheorem{prop}{Proposition}[section]
\newtheorem{lemma}[prop]{Lemma}
\newtheorem{thm}[prop]{Theorem}

\newtheorem{conj}[prop]{Question}
\newtheorem{question}[prop]{Question}
\newtheorem{coro}[prop]{Corollary} 

\newtheorem{rema}[prop]{Observation}
\newtheorem{defi}[prop]{Definition}
\newtheorem{example}[prop]{Example}
\newtheorem{fact}[prop]{Fact}

\newtheorem{theorem}{Theorem}

\newtheorem{definition}[theorem]{Definition}

\newcommand{\vertiii}[1]{{\left\vert\kern-0.25ex\left\vert\kern-0.25ex\left\vert #1 
    \right\vert\kern-0.25ex\right\vert\kern-0.25ex\right\vert}}

\begin{document}
    
\thanks{Research of W. Corrêa supported by São Paulo Research Foundation (FAPESP), Grants 2016/25574-8, 2018/03765-1, 2019/09205-0, 2021/13401-0, and by National Council for Scientific and Technological Development - CNPq - Brazil, Grant 304990/2023-0. Research of V. Ferenczi supported by São Paulo Research Foundation (FAPESP) Grants 2016/25574-8 and 2022/04745-0, and by National Council for Scientific and Technological Development (CNPq), Grants 303731/2019-2 and 304194/2023-9. Research of R. Gesing partially supported by São Paulo Research Foundation (FAPESP), Grants 2018/08362-2, 2017/14848-2, as well as by the Deutsche Forschungsgemeinschaft (DFG, German Research Foundation) under Germany’s Excellence Strategy – EXC 2044 – 390685587, Mathematics Münster – Dynamics – Geometry – Structure; the Deutsche Forschungsgemeinschaft (DFG, German Research Foundation) – Project-ID 427320536 – SFB 1442, and ERC Advanced Grant 834267 - AMAREC. Research of P. Tradacete partially supported by grants PID2020-116398GB-I00 and CEX2019-000904-S funded by MCIN/AEI/10.13039/501100011033, as well as by a 2022 Leonardo Grant for Researchers and Cultural Creators, BBVA Foundation.}

\author[W. Corr\^ea]{Willian Corr\^ea}
\address{Willian Corr\^ea, 
Departamento de Matem\'atica, Instituto de Ciências Matemáticas e de Computação, Universidade de São Paulo, Avenida Trabalhador São-carlense, 400 - Centro,
CEP: 13566-590, São Carlos, SP, Brazil.}
\email{willhans@icmc.usp.br}

\author[V. Ferenczi]{Valentin Ferenczi}
\address{Valentin Ferenczi, 
Departamento de Matem\'atica, Instituto de Matem\'atica e Estat\'istica, Universidade de S\~ao
Paulo, rua do Mat\~ao 1010, 05508-090 S\~ao Paulo SP, Brazil, and
Equipe d’Analyse Fonctionnelle, Institut de Math\'ematiques de Jussieu, Universit\'e Pierre et
Marie Curie - Paris 6, Case 247, 4 place Jussieu, 75252 Paris Cedex 05, France.}
\email{ferenczi@ime.usp.br}

\author[R. Gesing]{Rafaela Gesing}
\address{Rafaela Gesing, 
Mathematisches Institut,
University of M{\"u}nster, 
Einsteinstr.\ 62, 
48149 M{\"u}nster, Germany}
\email{rgesing@uni-muenster.de}

\author[P. Tradacete]{Pedro Tradacete}
\address{Pedro Tradacete, 
Instituto de Ciencias Matem\'aticas (CSIC-UAM-UC3M-UCM)\\
Consejo Superior de Investigaciones Cient\'ificas\\
C/ Nicol\'as Cabrera, 13--15, Campus de Cantoblanco UAM\\
28049 Madrid, Spain.}
\email{pedro.tradacete@icmat.es}

\title{Extremes of interpolation scales of Banach spaces}




\maketitle

\begin{abstract}
M. Daher gave conditions so that the spheres of the spaces in the interior of a complex interpolation scale  are uniformly homeomorphic. We look for sufficient conditions for the validity of this result and related ones on the extremes of the scale, with applications to uniform homeomorphism between spheres of Banach spaces and the sphere of $\ell_2$.
\end{abstract}

\tableofcontents

\section{Introduction}

\subsection{Uniform homeomorphisms of spheres}

The main open question we want to address in this paper is whether the unit sphere $S_X$ of any separable uniformly convex Banach space $X$ must be uniformly homeomorphic to the unit sphere of the (separable) Hilbert space.

This question is known to have a positive answer  when $X$ is assumed to have an unconditional basis.  Odell-Schlumprecht, \cite{OS}, proved that the sphere $S_X$ of $X$ is uniformly homeomorphic to the sphere of $\ell_2$ if 
$X$ has an unconditional basis and non-trivial cotype. This was extended by Chaatit to Banach lattices with a weak unit and with non-trivial cotype \cite{Cha}, and by Kalton to subspaces of superreflexive Banach lattices, see \cite[Corollary 9.11]{BL}.

Daher \cite{Da} and Kalton observed independently that complex interpolation could also be used in this context: i) if $X$ is uniformly convex and has an unconditional basis, then there exists a compatible couple of interpolation spaces $(\ell_2,Y)$, with $Y$ uniformly convex, and 
such that $X=(\ell_2,Y)_\theta$ for some $0<\theta<1$; ii) in this case generalizations of the Mazur maps witness that $S_X$ is uniformly homeomorphic to $S_{\ell_2}$.

The question remains widely open for spaces without lattice structure. In this case, we shall rely on an application of an interesting result of Daher \cite{D}.  Namely if $\ell_2$ and $X$ are both in the interior of an interpolation scale of uniformly convex spaces, without requirement of lattice structure, then the unit sphere of $X$ is uniformly homeomorphic to the unit sphere of $\ell_2$. Good illustrations of this are provided by the scales of $\ell_p$ spaces, Lebesgue spaces $L_p$ or $H_p$ spaces \cite{JJ}. 
Other examples of scales with $\ell_2$ in the interior of the domain are those associated to the pair $(\overline{X}^*,X)$, where
 $X$ has a Schauder basis and is uniformly convex \cite{Watbled} (here and throughout $\overline{X}$ denotes the complex conjugate of a complex Banach space $X$). 
 
Therefore, a first idea is extrapolating the scale $X_t$ generated by $(\overline{X}^*,X)$ to the right to make $X$ in the interior of a bigger scale which also contains $\ell_2$.
Since this seems to be a very difficult problem, a slightly less ambitious goal would be to find conditions so that the uniform homeomorphisms provided by Daher between $\ell_2=X_{1/2}$ and $X_t$, $t<1$, could be carried over to $X_1=X$. It is conceivable that these conditions do not involve lattice structure, since Daher's result does not rely on such properties.

This is why the main part of our paper is dedicated to the behaviour ``at the border'' of Daher's techniques.

\subsection{Spaces with Property (H)}

The same interpolation techniques could be relevant to another property of Banach spaces.
In \cite{KY}, Kasparov and Yu proved the Strong Novikov conjecture (with coefficients) for any group coarsely embeddable into a Banach space satisfying a geometric condition called Property (H). 

\begin{definition}[Kasparov-Yu]
    A real Banach space $X$ is said to have Property (H) if there exists an increasing sequence $(V_n)$ of finite dimensional subspaces of $X$ and an increasing sequence $(W_n)$ of finite dimensional subspaces of $\ell_2$ satisfying
        \begin{enumerate}
            \item $V = \bigcup\limits_{n=1}^{\infty} V_n$ is dense in $X$;
            \item If $W = \bigcup\limits_{n=1}^{\infty} W_n$ then there is a uniformly continuous map $\psi : S_V \rightarrow S_W$ whose restriction to each $S_{V_n}$ is a homeomorphism onto $S_{W_n}$.
        \end{enumerate}
\end{definition}

\begin{theorem}[Kasparov-Yu] Any group which is coarsely embeddable into a Banach space with Property (H) satisfies the Strong Novikov conjecture (with coefficients).\end{theorem}

Kasparov and Yu observed that uniformly convex spaces with (certain) unconditional bases have Property (H), as a consequence of results on the uniform homeomorphism classification of spheres of Banach spaces mentioned in Subsection 1.1. However, this observation does not really rely on unconditionality. As a first example, assume an interpolation scale of uniformly convex spaces with common $1$-Schauder basis contains 
both $X$ and $\ell_2$ in the interior of the scale. Then standard arguments imply that the Daher uniform homeomorphism from $S_X$ to $S_{\ell_2}$ preserves each finite dimensional subspace generated by an initial segment of the basis, and therefore witnesses that $X$ has property (H).
So again, in particular an appropriate extrapolation result for a space $X$ would guarantee that $X$ has Property (H). And since existence of extrapolation seems a lot to ask for, we are again led, in this case, to the study of Daher's techniques ``at the border'' and its continuity properties.
A possible advantage to take into account in this approach is that property (H) does not require uniform homeomorphism of the spheres; only continuity is required in one direction.

\subsection{Outline of the paper}

Given $(X_0, X_1)$ a compatible couple of Banach spaces, and $\theta \in (0,1)$, the complex interpolation space $X_\theta$ can be obtained as a quotient of a space $\mathcal F^{\infty}$ of analytic functions from the unit strip $\mathbb S=\{z \in \mathbb{C} : 0 \leq Re(z) \leq 1\}$ to $X_0 + X_1$, which satisfy certain integrability and measurability conditions. Daher shows in \cite{Da} that under the hypothesis of uniform convexity, the infimum in the quotient norm of an element in the unit sphere of the interpolation space $X_\theta$ is reached by a unique function in $\mathcal{F}^{\infty}$, which we shall call minimal function  (Definition \ref{def:optimal}). Minimal functions are a key element in Daher’s techniques to obtain uniform homeomorphisms between unit spheres of interpolation spaces. In Section \ref{sec:complexinterpolation} we recall the complex interpolation method and lay down the work of Daher and our setting to study minimal functions.

In Section \ref{sec:extrapolation} we introduce the notion of extrapolation of an interpolation scale of Banach spaces (Definition \ref{def:extrapolation}), a condition that when satisfied trivially implies the sought-after homeomorphism between unit spheres on the boundary of the scale. We also study the behaviour of Omega operators in interpolation theory and its connections with extrapolation. 

In Section \ref{sec:continuity} we introduce two notions of continuity for compatible Banach couples, that are formally weaker than extrapolation, namely, continuity to the right/left (Definition  \ref{def:continuos_to_the_right}) and right/left norm continuity (Definition \ref{def:right_norm_continuous}). These concepts provide insightful information on the minimal functions and Omega operators.

Section \ref{sec:vertical} is dedicated to the study of the vertical symmetries intrinsic to the complex interpolation method. By studying the periodic behaviour of the so-called vertical orbits, we deduce that periodicity imposes heavy restrictions on minimal functions (Proposition \ref{cndn}). In addition, we define the limit Mazur map (Definition \ref{def:mazurmap}), a candidate for being a uniform homeomorphism between unit spheres, and a map witnessing property (H). We note that its behaviour depends on the vertical orbits. 

Relying on the use of the vertical homeomorphisms considered in the previous section, in Section \ref{sec:uniform_scales} we introduce the notion of {uniformities} (Definition \ref{def:uniformities}), which are properties of interpolation scales concerning the modulus uniform continuity of these vertical mappings, which in turn ensure favorable properties of the limit Mazur map. Through uniformity within the scale, we obtain uniform homeomorphism between the spheres of the spaces within the domain's interior and its boundary:

\begin{theorem}[Corollary \ref{cor:main_result}]
  Let $X$ be a uniformly convex complex space with a $1$-FDD. 
Consider the interpolation scale $(\overline{X}^* ,X)$. If this scale is  uniform to the right, then $S_X$ is uniformly homeomorphic to $S_{\ell_2}$.
\end{theorem}

\section{Background on complex interpolation and Daher's Theorem}\label{sec:complexinterpolation}

Recall that a couple $(X_0, X_1)$ of Banach spaces is said to be \emph{compatible} if there are a Banach space $U$ and continuous linear injections from $X_0$ and $X_1$ into $U$. In practice, we substitute $U$ by the sum space $X_0 + X_1$, and so we may consider the injections as continuous inclusions.

We now give the definition of the spaces of analytic functions we use for interpolation. These are direct generalizations of the spaces defined by Daher in \cite{Da}.

Let $\mathbb{S} = \{z \in \mathbb{C} : 0 \leq Re(z) \leq 1\}$, and for $z \in \overset{\circ}{\mathbb{S}}$ (the interior of the strip $\mathbb{S} $) let $\mu_{z}$ be the Poisson kernel for $\mathbb{S}$ with respect to $z$ (cf. \cite{Widder}). Also, for $j = 0, 1$ let $\mu_{z}^j$ be measure on $\R$ defined by $\mu_{z}^j(A)=\mu_{z}(j+iA)$, i.e. $\mu_z^j$ is the ``restriction" of $\mu_z$ to the vertical line of real part equal to $j$.

For $1<p<+\infty$, and $z \in \overset{\circ}{\mathbb{S}}$, let $\mathcal{F}_{z}^p$ be the space of functions $f : \mathbb{S} \rightarrow X_0 + X_1$ such that
\begin{enumerate}
    \item $f$ is analytic on the interior of $\mathbb{S}$;
    \item The functions $t \mapsto f(j + it)$, $j = 0, 1$, are (Bochner) measurable from $\mathbb{R}$ into $X_j$;
    \item
        \[
        \|f\|_{\mathcal{F}_{z}^p} := \Big(\int\limits_{-\infty}^{\infty} \|f(it)\|_{X_0}^p d\mu_{z}^{0}(t) + \int\limits_{-\infty}^{\infty} \|f(1 + it)\|_{X_1}^p d\mu_{z}^{1}(t)\Big)^{\frac{1}{p}} < \infty
        \]
    \item $f$ satisfies the Poisson kernel representation
        \[
        f(w) = \int_{\partial \mathbb{S}} f(w') d\mu_{w}(w')
        \]
        for every $w$ in the interior of $\mathbb{S}$.
\end{enumerate}

For each $z \in \overset{\circ}{\mathbb{S}}$, the complex interpolation space is defined as $X_{z} = \{f(z) : f \in \mathcal{F}_{z}^p\}$, endowed with the quotient norm, which we denote by $\Vert \cdot \Vert_z$. One can check that $X_z = X_{Re(z)}$ isometrically, so in general we restrict ourselves to considering $X_{\theta}$ for $\theta \in (0, 1)$. It can also be proved that this construction (for any choice of $p$ in $[1,+\infty)$) gives the same interpolation space as Calder\'on's complex interpolation method up to an equivalent norm (see \cite{Da}).

The space $\Finf$ is defined in the same way, replacing (3) by 
 \[
        \|f\|_{\Finf}:= \max(\|f(i(\cdot))\|_{L^{\infty}(X_0)}, \|f(1 + i(\cdot))\|_{L^{\infty}(X_1)}) < \infty.
        \]

Note that for $z_1, z_2 \in \overset{\circ}{\mathbb{S}}$ we have $\mathcal{F}_{z_1}^p = \mathcal{F}_{z_2}^p$ with equivalence of norms. 

\begin{defi} A compatible couple $(X_0, X_1)$ of Banach spaces is {\em regular} if the intersection $X_0 \cap X_1$ is dense in both $X_0$ and $X_1$.
\end{defi}

Regularity is an important property of compatible couples as shown by Daher's Theorem:

\begin{thm}\label{thm:Daher}[\cite{Da}, Proposition 3]
Let $(X_0, X_1)$ be a regular compatible couple of reflexive Banach spaces.
Let $\theta \in (0, 1)$ and $x \in S_{X_{\theta}}$.
Then
\begin{itemize}
    \item[(1)] there exists $G \in \Finf$ such that $G(\theta) = x$ and $\|G(z)\|_{z}=1$ for almost all $z \in \partial \mathbb{S}$, and so $\|G\|_{\mathcal{F}^\infty} = 1$;
    \item[(2)] if $X_0$ or $X_1$ is strictly convex then such $G$ is the unique $ G \in \mathcal{F}_{\theta}^2$ such that $G(\theta) = x$ and $\|G\|_{\mathcal{F}_{\theta}^2} = 1$. We denote such $G$ by $F_x^{\theta}$;
    \item[(3)] if $X_0$ or $X_1$ is uniformly convex, and $\theta, \theta'$ are interior to the scale, then the map $\phi_{\theta, \theta'}: x \mapsto F^{\theta}_x(\theta')$ is a uniform homeomorphism between $S_{X_{\theta}}$ and $S_{X_{\theta'}}$.
\end{itemize}

\end{thm}

Note that the previous theorem was generalized to other interpolation methods in \cite{Co}.

\begin{defi} \label{def:optimal}
 Given $(X_0, X_1)$ a regular compatible couple, $\theta \in (0,1) $ and $x \in S_{X_\theta}$, a map  $G \in \mathcal{F}^\infty$ such that $G(\theta) = x$ and $\|G(z)\|_z=\|x\|_\theta$ for almost all $z \in \partial \mathbb{S}$ will be called a {\em minimal function} for $x$ in $\theta$.
 
 A compatible couple admitting a unique minimal function for any $0<\theta<1$ and any $x \in X_\theta$ will be called {\em optimal} \cite{Stability}.
\end{defi}

In practical cases and to study the behaviour of the scale close to the right side of $\partial \mathbb{S}$, it is  usually not a loss of generality, and often useful, to assume that $X_0$ and $X_1$ are both uniformly convex. Indeed under this assumption $X_\theta$ is also uniformly convex for all $0<\theta<1$ and we may just look at the scale generated by $(X_\theta,X_1)$ instead of $(X_0,X_1)$.

\section{Extrapolation}\label{sec:extrapolation} 
{}
\subsection{Scales admitting extrapolation}
The uniform homeomorphism between spheres obtained by Daher holds for spaces in the interior of the scale (i.e. for $\theta, \theta' \in (0,1)$ in Theorem \ref{thm:Daher}). 
We investigate a first condition which would allow to extend it to the extremes of the scale. Therefore we introduce the following definition:

\begin{defi} \label{def:extrapolation} We say that
a regular compatible couple $(X_0,X_1)$ admits {\em extrapolation on the right} (\emph{resp. on the left}) if there exist some space $Y$ such that $(X_0,Y)$ (resp. $(Y, X_1)$) is a compatible couple, and some $\theta$ such that $(X_0,Y)_\theta=X_1$ (resp. $(Y, X_1)_{\theta} = X_0$) isometrically. 
\end{defi}

Note that we may always assume that $(Y, X_1)$ or $(X_0, Y)$ is a regular compatible couple, since $(Y, X_1)_{\theta} = (\overline{Y \cap X_1}^{Y}, X_1)_{\theta}$ and $(X_0, Y)_{\theta} = (X_0, \overline{Y \cap X_0}^{X_0})_{\theta}$.
Note also that if $(X_0, X_1)$ admits extrapolation on the right and  $X_1$ are uniformly convex, then $Y$ may be assumed to be uniformly convex as well (if it is not then replace it by $(X_1,Y)_{1/2}$ up to changing the value of $\theta$, \cite{CR}). In particular Daher's result implies the following:

\begin{prop} Assume $(X_0,X_1)$ is a regular compatible couple admitting extrapolation to the right. If $X_1$ is uniformly convex then the map $x \mapsto F_x^{\theta}(1)$ defines a uniform homeomorphism between the unit spheres of $X_\theta$
and $X_1$.\end{prop}

Kalton \cite{KaltonDifferentials, kaltonhardy} proved that scales of uniformly convex K\"othe function spaces admit extrapolations up to renorming. Let us state this result here:

\begin{prop}\label{thm:Kotheextrapolation1}{\cite[Proposition 2.1]{kaltonhardy}}
Let $(X_0, X_1)$ be a compatible couple of uniformly convex K\"othe function spaces. Then there is $\theta \in (0, 1)$ and a K\"othe function space $Y$ such that $(X_0, Y)_{\theta} = X_1$ up to equivalence of norms.
\end{prop}

Note that in the previous result it is enough to suppose that $X_1$ is uniformly convex, since in this case $X_{\theta}$ is uniformly convex for every $\theta \in (0, 1)$. Therefore, using Daher's results we deduce the following:

\begin{coro}
Let $(X_0,X_1)$ be a compatible couple of K\"othe function spaces such that $X_1$ is uniformly convex. Then the unit spheres of $X_1$ and $X_\theta$ are uniformly homeomorphic for every $\theta \in (0, 1)$. 
\end{coro}

The uniform homeomorphism is obtained through interpolation of renormings of $X_0$ and $X_1$ (noting that equivalent renormings preserve uniform homeomorphisms of unit spheres).

\medskip

It is an open question whether this extends to  uniformly convex Banach spaces with weaker structure, for example:

\begin{question}\label{question:extrapolation} Let $(X_0,X_1)$ be a regular pair of spaces with a common $1$-Schauder basis, or more generally with a common finite dimensional decomposition with constant $1$ (1-FDD for short, cf. \cite[1.g]{LT}). If $X_1$ is uniformly convex then does (some renorming of) the pair admit extrapolation on the right?
\end{question}

Note that if Question \ref{question:extrapolation} has an affirmative answer then we would obtain for any uniformly convex space $X$ with a 1-FDD:
\begin{itemize}
 \item Uniform homeomorphism between $S_X$ and $S_{\ell_2}$.
 \item Property (H) for $X$.
\end{itemize}
 To see this just consider the compatible couple
 $(\overline{X}^*,X)$ and use the fact that
 $(\overline{X}^*,X)_{1/2}=\ell_2$.
 Therefore the Novikov conjecture would hold for all superreflexive spaces with a 1-FDD.
 Note that it is an open question whether every separable uniformly convex space embeds into one with an FDD (or even a basis) \cite[Problem 2.16]{OS2}.

Note also that, quite surprisingly, extrapolation does not have to be unique: there exist different complex interpolation scales that agree on an open subinterval of $[0, 1]$. See \cite{JNPZ, Kellogg} and \cite[pages 320-321]{Triebel}. In the case of scales of K\"othe function spaces that cannot happen, due to the results of \cite{KaltonDifferentials}.

\subsection{Extrapolation and Omega operators}

\medskip

Let us recall the definition of the $\Omega$ operators in interpolation theory. We have many choices of spaces of analytic functions to work with, so let us call any one of them $\mathcal{C}.$ Given a compatible couple $(X_0, X_1)$, $C >0$ and $\theta \in (0, 1)$, a homogenous map $B_{\theta} : X_{\theta} \rightarrow \mathcal{C}$ is called a \emph{C-selector} if for every $x \in X_{\theta}$:
\begin{enumerate}
    \item $\|B_{\theta}(x)\|_{\mathcal{C}} \leq C\|x\|_{\theta}$;
    \item $B_{\theta}(x)(\theta) = x$.
\end{enumerate}
For most purposes we can suppose that $B_{\theta}$ is defined only on a dense subspace of $X_0 \cap X_1$. For example, for the couple $(\ell_{\infty}, \ell_1)$ we have a $1$-selector which on vectors of finite support with $\|x\|_{\ell_{\frac{1}{\theta}}} = 1$ is given by
\[
B_{\theta}(x)(z) = \sum sgn(x_n) \abs{x_n}^{z/\theta} e_n.
\]
Whenever it is possible to take a $1$-selector, we do so.

Given a $C$-selector $B_{\theta} : X_{\theta} \rightarrow \mathcal{C}$, for each $n \geq 1$ we let $\Omega^{(n)}_{\theta}(x) = \frac{B_{\theta}(x)^{(n)}(\theta)}{n!}$, i.e., $\Omega^{(n)}_{\theta}(x)$ is the $n$-th Taylor coefficient of $B_{\theta}(x)$. A particular notation is $\Omega_{\theta}(x) = \Omega^{(1)}_{\theta}(x)$.

For the case of the couple $(L_{\infty}(\mu), L_1(\mu))$, the $\Omega$ operators are the Kalton-Peck maps on $L_p(\mu)$, $1<p<\infty$, for which
$$\Omega_\theta(f)= \frac{1}{\theta} f \cdot log(|f|/\|f\|),$$
with $\theta=\frac{1}{p}$, \cite{KP}. 

The $\Omega$ operators are related to the theory of twisted sums and commutator estimates, and in the case of superreflexive K\"othe function spaces, Kalton's \cite[Theorem 7.6]{KaltonDifferentials} ensures that for any given $\theta$ the space $X_{\theta}$ and the map $\Omega_{\theta}$ completely determine the spaces $X_0$ and $X_1$.

Since the $\Omega$ operators are the Taylor coefficients of optimal functions, one strategy would be to study their growth to try to obtain that optimal functions may be extended beyond the strip.

For $n \in \mathbb{N}$ and $t \in [0, 1]$ let $X_{n, t} = (\mathbb{C}, e^{-n t} \left|\cdot\right|)$ ($\mathbb{C}$ with a weight). Then $(X_{n, t})_{t \in [0, 1]}$ is a complex interpolation scale, and $\Omega_{t}(x) = n x$. Note that all the spaces have the same modulus of uniform convexity, and this family of examples shows that we cannot expect a bound on the norm of $\Omega_t(x)$ in terms of the norm of $x$ and the modulus of uniform convexity.
Still, we have extrapolation of those scales to the whole real line! Maybe the following question has a satisfactory answer (it covers the families $(X_{n, t})$):

\begin{question}
Suppose that $\Omega_{\theta} : X_{\theta} \rightarrow X_{\theta}$ is bounded (possibly uniformly bounded on an interval). Under what conditions do we get extrapolation?
\end{question}

A possible strategy would be the following: prove that if $\|\Omega_{\theta}(x)\|_{\theta} \leq C$ then $\|B_{\theta}(x)^{(n)}(\theta)\|_{\theta} \leq C^n$. Then $\Big(\frac{\|B_{\theta}(x)^{(n)}(\theta)\|_{\theta}}{n!}\Big)^{\frac{1}{n}} \leq \frac{C}{\sqrt[n]{n}} \rightarrow 0$. Therefore, the radius of convergence of the Taylor series would be infinite, and we can extend the optimal function.

\

Note that if $X$ is a space for which $(\overline{X}^*,X)_{1/2}=\ell_2$, and if this interpolation scale admits extrapolation to the right, then $X$ has to be $\theta$-Hilbertian (see the definition in \cite{Pisier}). Pisier uses those spaces to characterize the interpolation space $(B_r(L_2(\mu)), B(L_2(\mu))^{\theta}$, where $B_r(L_2(\mu))$ is the space of regular operators on $L_2(\mu)$. He also gives a characterization of the uniform curved spaces for which the corresponding modulus $\Delta_X(\varepsilon)$ is $O(\varepsilon^{\alpha})$ for some $\alpha > 0$. It is still unknown whether all uniformly convex spaces are $\theta$-Hilbertian. Therefore, Question \ref{question:extrapolation} in full generality seems very difficult, and this is why we are looking for properties which are weaker than admitting extrapolation.
One first family of such properties are {\em continuities}.

\section{Continuity of interpolation scales}\label{sec:continuity}

\subsection{A first notion of continuity of scales}

\begin{defi}  \label{def:continuos_to_the_right}  Let $(X_0, X_1)$ be a regular compatible couple, and let $X \subset X_0 \cap X_1$ be a dense subspace. We will say that the scale is \emph{continuous to the right} (resp. \emph{to the left}) in $X$ if every $x \in X$ admits a unique minimal function $F_x^{\theta}$ at $\theta$ for any $0 < \theta < 1$, and that function is continuous on the semiclosed strip
 $\{z: 0<Re(z) \leq 1\}$ (resp. $\{z: 0 \leq Re(z) < 1\}$)
 as a map with values in the normed space $X_0+X_1$.
\end{defi}

 When $(X_0,X_1)$ is a scale of spaces with a common 1-basis or, more generally a common 1-FDD, let $(Y_n)_n$ denote the sequence of common finite dimensional spaces in the decomposition. We shall assume in this case that $X$ is $c_{00}(Y_n)$, the linear span of $\bigcup_n Y_n$. Thus, continuity in this case means that every $x$ finitely supported with respect to the decomposition $(Y_n)$ admits a unique minimal function $F_x^\theta$, and that this function is continuous with values in $X_0+X_1$.

It follows from the definition of minimal functions that if a scale of finite dimensional spaces is continuous to the right then $F_x^\theta(z)$ takes values in $S_{X_1}$ whenever $Re(z)=1$ (see Proposition \ref{props:Delta_theta_in_S_theta}). The same holds for a scale with a common 1-FDD, for all $x \in X$ as above.

An interesting example is the scale of $L_p$ spaces on $[0, 1]$ when the endpoint to the right is $L_1$. Let us fix the left endpoint as $L_{p_0}$ for some $1 < p_0 < \infty$. Then the pair $(L_{p_0}, L_1)$ is regular. Define $p(z)$ on the closed strip by $\frac{1}{p(z)} = \frac{1 - z}{p_0} + z$. Then $(L_{p_0}, L_1)_{\theta} = L_{p(\theta)}$ and for $x \in L_{p(\theta)}$ of norm $1$ in a dense subspace of $L_{p_0} \cap L_1$ we have
\[
F_x^{\theta}(z) = \left|x\right|^{\frac{p(\theta)}{p(z)}} \frac{x}{\left|x\right|}
\]

Note that in this case the Mazur map is a uniform homeomorphism between the spheres. 
However $L_1$ is not uniformly convex or even reflexive.

For most natural scales, continuity is formally weaker than extrapolation:

\begin{lemma} If a scale of spaces with a common 1-FDD and with either $X_0$ reflexive strictly convex or $X_1$ uniformly convex admits extrapolation to the right, then it is continuous to the right.
\end{lemma}
\pf
If $X_0$ is reflexive strictly convex and $(X_0, X_1)$ admits extrapolation to the right then by \cite{BeLo} pp 105, note 4.9.1-4.9.6, we have that $X_1$ is reflexive. Since both $X_0$ and $X_1$ are reflexive and one is strictly convex, Daher's result applies to obtain uniqueness of the minimal function. Continuity along the vertical line corresponding to $X_1$ follows from analyticity.

If $X_1$ is uniformly convex then $X_\eta$ is uniformly convex and therefore reflexive strictly convex for any $\eta>0$; therefore the previous case applies to the scales $\eta < Re(z) < 1$ and therefore to the whole scale.
\pff

\

We have the following criteria:

\begin{prop} Suppose we have a compatible couple $(X_0, X_1)$ with a common 1-FDD and $x \in X_0 \cap X_1$ admits a unique minimal function $F_{x}^{\theta}$ which has image in a fixed finite dimensional space $X$ contained in $X_0 \cap X_1$ for every $\theta$. 
Consider the application $\Gamma_x : \theta \mapsto F_{x}^{ \theta}$, and suppose that $\Gamma_x$ is continuous from $(0, 1)$ into $\mathcal{F}^{\infty}$. Then the function 
$\theta \mapsto \|\Omega_{\theta}(x)\|_{\theta}$ is continuous.

In particular, if in addition we have extrapolation to the right, then
\begin{itemize}
    \item[(1)]
$$\sup\limits_{\frac{1}{2} < \theta < 1} \|\Omega_{\theta}(x)\|_{\theta} < \infty$$
\item[(2)]
$$\sup\limits_{\frac{1}{2} < \theta < 1} \left|\left(\frac{d}{dt}\|x\|_{t}\right)\Big|_{{t = \theta^-}}\right| < \infty$$
\end{itemize}
\end{prop}

\pf
 We can do the proof under slight more generality. Note first we have continuity of the following maps:
 \[
 \Gamma_x \times Id : (0, 1) \rightarrow \mathcal{F}^{\infty} \times (0, 1)
 \]
 given by $\theta \mapsto (F_x^{\theta}, \theta)$, and
 \[
 \delta^{\prime} : Y \times (0, 1) \rightarrow X \times (0,1)
 \]
 given by $(f, \theta) \mapsto (f^{\prime}(\theta), \theta)$. Here $Y$ is the linear span of $\{F_x^{\theta} : \theta \in (0, 1)\}$ and $X$ is the finite dimensional space where all $F_x^{\theta}$ have their image. Let us explain why $\delta^{\prime}$ is continuous.
 
 Let $(f_1, \theta_1), (f_2, \theta_2) \in Y \times (0, 1)$. Continuity of $\delta^{\prime}$ boils down to estimating
 \[
 \|f_1^{\prime}(\theta_1) - f_2^{\prime}(\theta_2)\|_X \leq \|f_1^{\prime}(\theta_1) - f_2^{\prime}(\theta_1)\|_X + \|f_2^{\prime}(\theta_1) - f_2^{\prime}(\theta_2)\|_X
 \]
 The evaluation of the derivative $\delta_{\theta}^{\prime} : \mathcal{F}^{\infty} \rightarrow X_0 + X_1$ is bounded:
 \[
 \|f^{\prime}(\theta)\|_X = \frac{1}{2\pi} \norm{ \int_C \frac{f(z)}{(z - \theta)^2} dz}_X \leq K \|f\|_{\infty}
 \]
 where $C$ is some simple closed curve around $\theta$ and $K$ is a constant that depends on $\theta$ and our choice of $C$. Therefore, the first term of the right side is taken care of. For the second, we consider another simple closed curve, this time around $\theta_1$ and $\theta_2$, and we have:
 \begin{eqnarray*}
 \|f^{\prime}(\theta_1) - f^{\prime}(\theta_2)\|_X & = & \frac{1}{2\pi} \norm{\int_C \frac{f(z)}{(z-\theta_1)^2} - \frac{f(z)}{(z-\theta_2)^2} dz}_X \\
 & = & \frac{1}{2\pi} \norm{ \int_C \frac{(z-\theta_2)^2 - (z-\theta_1)^2}{(z-\theta_1)^2(z-\theta_2)^2} f(z) dz }_X \\
 & \leq & K \sup_{z \in C} \abs{\frac{(z-\theta_2)^2 - (z-\theta_1)^2}{(z-\theta_1)^2(z-\theta_2)^2}} \|f\|_{\infty}
 \end{eqnarray*}
 and therefore $\|\delta_{\theta_1}^{\prime} - \delta_{\theta_2}^{\prime}\| \rightarrow 0$ as $\theta_2 \rightarrow \theta_1$. This takes care of the second term of the inequality, and therefore $\delta^{\prime}$ is continuous.

 Now consider
 \[
 g : X \times (0, 1) \rightarrow \mathbb{R}
 \]
 given by $g(y, \theta) = \|y\|_{\theta}$. Let us prove that $g$ is continuous.
 
 Suppose $y_n \rightarrow y$ and $\theta_n \rightarrow \theta$. We have
 \[
 \|y\|_{\theta} - \|y_n\|_{\theta_n} = (\|y\|_{\theta} - \|y_n\|_{\theta}) + (\|y_n\|_{\theta} - \|y_n\|_{\theta_n})
 \]
 The first part goes to zero because we are in a finite dimensional space, and $\|\cdot\|_{\theta}$ is equivalent to the norm of $X$. The second part is trickier: using \cite[Theorem 3.1]{RochbergFamilies} and reiteration, we are able to find a constant $C > 0$ such that for $\eta$ close enough to $\theta$ we have for any $w \in X$
 \[
 \left|\|w\|_{\theta} - \|w\|_{\eta}\right| \leq C \left|\varphi(\eta)\right| \|w\|_{X_0 \cap X_1}
 \]
 where $\varphi$ is a conformal map from the strip to the unit disk of radius 1.1 sending $\theta$ to $0$. Indeed, let $\D_{1.1}$ be the open disk of radius 1.1. \cite[Theorem 3.1]{RochbergFamilies} tells us that if we have an interpolation family $(B_z)_{z \in \D_{1.1}}$ of finite dimensional spaces of same dimension, then there is a constant $c > 0$ such that
 \[
 \abs{\|w\|_0 - \|w\|_z} \leq c \abs{z} \int_{\mathbb{T}} \|w\|_{e^{i\theta}} d\theta
 \]
for every $w$. The norms on $X$ induced by $(X_0, X_1)_{\theta}$ and by $(X^{\|\cdot\|_0}, X^{\|\cdot\|_1})_{\theta}$ are equal for every $\theta \in (0, 1)$, because of the 1-FDD. The claim follows from considering the interpolation family $(X^{\|\cdot\|_{\varphi^{-1}(z)}})_{z \in \mathbb{D}_{1.1}}$.
 
Since $y_n \rightarrow y$ and we are in a finite dimensional space, $(\|y_n\|_{X_0 \cap X_1})$ is bounded, and we must have that $\|y_n\|_{\theta} - \|y_n\|_{\theta_n} \rightarrow 0$.

Now $\theta \mapsto \|\Omega_{\theta}(x)\|_\theta$ is the composition of those maps, and then it must be bounded on any compact subset of $(0, 1)$. 
\medskip
  
(1) follows if in addition we have extrapolation. 
\medskip
  
(2) is a consequence of \cite[Lemma 5.5]{Stability}.
\pff

\

The following example shows that we can have continuity without extrapolation.

\begin{example}
Consider the couple $(\ell_{\infty}, \ell_1)$, and let $x = e_1 + e_2$. Then
\[
F_x^{\theta}(z) = \|x\|_p \Big( \frac{1}{\|x\|_p^{\frac{z}{\theta}}} e_1 + \frac{1}{\|x\|_p^{\frac{z}{\theta}}} e_2 \Big)
\]
is the minimal function of $x$ at $\theta = \frac{1}{p}$, and
\[
\Omega_{\theta}(x) = \log \|x\|_{\theta}^{-1} x.
\]
 Then $F_x^{\theta}$ is continuous on the left and on the right line.
We have $$\lim_{\theta \rightarrow 0^+} \|\Omega_{\theta}(x)\|_{\theta} = \infty,$$ so there is no extrapolation to the left. On the other hand, $$\lim_{\theta \rightarrow 1^+} \|\Omega_{\theta}(x)\|_{\theta} < \infty,$$ and we have extrapolation to the right if we allow  the quasi-Banach spaces $\ell_p, 0<p<1$, in our definition (but not otherwise).
\end{example}
  
 \

Apart from this, what can be said about the relations between the different notions remains fairly open. For example:

\begin{question}

\

If a regular scale is continuous to the right, and $X_1$ is uniformly convex,  does this imply that the sphere of $X_1$ is uniformly homeomorphic to the sphere of $X_\theta$, $0<\theta<1$?

\

Can we find a regular scale of interpolation, with $X_1$ uniformly convex, which is  continuous to the right but which does not admit extrapolation to the right?

 \
 
 Is a regular scale of uniformly convex spaces necessarily continuous? 
\end{question}

\subsection{Norm continuous scales}

We now investigate another relevant notion of continuity.

We always have that the map $t \mapsto \|x\|_t$ is continuous for $t \in (0, 1)$ and $x \in X_0 \cap X_1$. The admittance of extrapolation implies the continuity of this map at the extreme points of the interval, as we observe below.

\begin{defi} For a regular compatible couple $(X_0,X_1)$ and for
 $x \in \Delta \coloneqq X_0 \cap X_1$, for $j=0,1$ let 
$$\vertiii{x}_j \coloneqq \lim_{t \rightarrow j} \|x\|_t.$$
\end{defi}

Note that the above limits exist by log-convexity of the map $t \mapsto \|x\|_t$ for all $x \in X_0 \cap X_1$ \cite[Lemma 5.1]{Stability}.

\begin{defi} \label{def:right_norm_continuous} 
For a regular compatible couple $(X_0,X_1)$ we say that the scale is \emph{right norm continuous} (resp. \emph{left norm continuous})  \emph{in $x \in \Delta\subset X_0 \cap X_1$}
if   $$\vertiii{x}_1=\|x\|_1 \;\;\;\;\;\; \text{(resp. } \vertiii{x}_0 = \|x\|_0 \text{)}$$ and that it is \emph{right norm continuous} (resp. \emph{left norm continuous}) if it is right (resp. left) norm continuous  in $x$ for all $x \in \Delta$. \end{defi}

\begin{rema} 
If a scale admits extrapolation on the right, then it is right norm continuous.
\end{rema}

\pf This is because we always have continuity in the interior of a scale \cite[Lemma 5.1]{Stability}. \pff

\

At first view, the above may seem like a criteria for the existence of extrapolation; however, we now show that all scales with a common 1-FDD admit norm continuity.

Given a compatible couple $(X_0,X_1)$ of Banach spaces, for $x\in X_0+X_1$ and $t>0$, we recall that the $K$-functional is defined by
$$
K(x,t;X_0,X_1)=\inf\{\|x_0\|_{X_0}+t\|x_1\|_{X_1}:x=x_0+x_1, x_i\in X_i\}.
$$
Following the terminology of Milman \cite{M} we say that:

\begin{defi} A  couple $(X_0,X_1)$ is normal in $x \in \Delta$ if 
\begin{enumerate}
\item $\lim_{t\rightarrow\infty} K(x,t;X_0,X_1)=\|x\|_{X_0}.$
\item $\lim_{t\rightarrow0} \frac{K(x,t;X_0,X_1)}{t}=\|x\|_{X_1},$
\end{enumerate}
\end{defi}

Recall that the Gagliardo completion of $X_0$ in $X_0+X_1$ is defined as $\tilde X_0=\{x\in X_0+X_1:\sup_t K(x,t;X_0,X_1)<\infty\}$, endowed with the norm $\|x\|_{\tilde X_0}=\sup_t K(x,t;X_0,X_1)$. Analogously, the Gagliardo completion of $X_1$ in $X_0+X_1$ is $\tilde X_1=\{x\in X_0+X_1:\sup_t \frac{K(x,t;X_0,X_1)}{t}<\infty\}$, endowed with the norm $\|x\|_{\tilde X_1}=\sup_t \frac{K(x,t;X_0,X_1)}{t}$.

According to \cite[Theorem 5.1.5]{BennetSharpley}, for $x\in X_0+X_1$ and $t>0$ we have
$$
K(x,t;X_0,X_1)=K(x,t;\tilde X_0,\tilde X_1).
$$
Hence, $(\tilde X_0,\tilde X_1)$ is normal in every $x \in \Delta$.

\begin{rema}[\cite{M}]\label{obs:Milman} If $(X_0,X_1)$ is normal  in $x \in \Delta$ then it is norm continuous in $x$. \end{rema}

The finite dimensional case in the following result also follows from \cite{Rochberg_Weiss_Analytic_Families}, (2.9), but we give a proof for the sake of completeness.

\begin{prop}\label{310} If $(X_0,X_1)$ is a regular couple which is finite-dimensional, or with a common $1$-FDD, then it is norm continuous. 
\end{prop}

\pf
Let $(Y_n)$ be a common 1-FDD. We prove that it is normal on every vector in $c_{00}(Y_n)$, therefore norm continuous on the linear span of the finite dimensional spaces in the decomposition, and therefore on $\Delta$.
Suppose $x\in c_{00}(Y_n)$, we claim that for $j=0,1$ $$\|x\|_{\tilde X_j}=\|x\|_{X_j}.$$
Indeed, for $x\in X_0$ and every $t>0$ we have $K(x,t;X_0,X_1)\leq \|x\|_{X_0}$. Thus, $\|x\|_{\tilde X_0}\leq \|x\|_{X_0}$. For the converse inequality, using \cite[Theorem 5.1.4]{BennetSharpley} we have
$$
\|x\|_{\tilde X_0}=\inf\{\sup_n\|y_n\|_{X_0}: (y_n) \subset X_0, \|x-y_n\|_{X_0+X_1}\rightarrow0\}.
$$

Suppose that $\max(\textrm{supp}(x))=m$ and let $P_m$ denote the corresponding FDD projection (on $X_0$). Given any $(y_n)\subset X_0$ such that $\|x-y_n\|_{X_0+X_1}\rightarrow 0$, we have that $\|x-P_my_n\|_{X_0+X_1}\rightarrow0$. Since these elements belong to certain finite dimensional subspace, we have $\|x- P_my_n\|_{X_0}\rightarrow 0$, so in particular $$\sup_n \|P_m y_n\|_{X_0}\geq\|x\|_{X_0}.$$ Since, $\|P_m\|\leq 1$ we get
$$
\|x\|_{X_0}\leq \sup_n\|y_n\|_{X_0},
$$ and taking infima on $(y_n)$ it follows that $\|x\|_{X_0}\leq \|x\|_{\tilde X_0}$. Similarly, one can show that $\|x\|_{\tilde X_1}=\|x\|_{X_1}$.

The conclusion follows from Milman's result (Observation \ref{obs:Milman}) and the fact that $(\tilde X_0,\tilde X_1)$ is normal in $\Delta$.
 \pff

\begin{question}
Does there exist a regular compatible couple with $X_1$ uniformly convex and failing norm continuity to the right?
\end{question}

\subsection{An observation regarding the limit}

\begin{lemma}\label{limit}
Let $(X_0,X_1)$ be a regular compatible couple. 
Let $F: \mathbb S \rightarrow X_0+X_1$ be analytic, such that
$F(s)$ belongs to the same finite-dimensional subspace for all $0<s<1$, and such that $F(1)=\lim_{s \rightarrow 1}F(s)$. Then $$\lim_{s \rightarrow 1}\|F(s)\|_s=\|F(1)\|_1.$$
\end{lemma}

\pf
Indeed 
$$ 
|\|F(s)\|_s - \|F(1)\|_s| \leq \|F(s)-F(1)\|_s \leq K\|F(s)-F(1)\|_{X_0+X_1}$$ for some fixed constant $K$, and this tends to $0$ by continuity of $F$.
Then use Proposition \ref{310} to see that
$|\|F(s)\|_s-\|F(1)\|_1| \leq |\|F(s)\|_s - \|F(1)\|_s|+\varepsilon $  for $s$ close enough to $1$. \pff

We deduce the following corollary.

\begin{coro}\label{corocoro}
Let $(X_0,X_1)$ be a regular compatible couple with a common $1$-FDD.  Let $0<\theta<1$ and $x,y$ in the linear span of the spaces inducing the FDD such that $F_x^{\theta}(s)$ and $F_y^{\theta}(s)$ admit limits when $s \rightarrow 1$.
Then
$$\lim_{s \rightarrow 1}\|F_x^\theta(s)-F_y^\theta(s)\|_s=\|F_x^\theta(1)-F_y^\theta(1)\|_1.$$
\end{coro}

\

\section{Vertical maps, uniform continuities, and limit maps}\label{sec:vertical}

We now carry on a more subtle study of some characteristics of Daher's method. Our objective is still to find some weaker conditions than extrapolation that imply some good limit behavior of the scale. It will be important to consider all points in the interior of the strip, and not only in the interval $(0,1)$, and to exploit the vertical symmetries of the construction.
We first recall Daher's theorem and fix some notation. 

\begin{thm}[\cite{Da}, Proposition 3]\label{daher}
Let $(X_0, X_1)$ be a regular compatible couple of reflexive Banach spaces. If $X_0$ or $X_1$ is strictly convex then for every $z$ in the interior of $\mathbb{S}$ and for every $x \in S_{X_{z}}$ there is a unique $F_x^{z} \in \mathcal{F}_{z}^2$ such that $F^{z}_x(z) = x$ and $\|F^{z}_x\|_{\mathcal{F}_{z}^2} = 1$. We have that
$\|F^z_x(it)\|_0=\|F^z_x(1+it)\|_1=1$ for a.e. $t \in \mathbb R$.

Furthermore, if $X_0$ or $X_1$ is uniformly convex, then the map $\Gamma_{z}: x \mapsto F^{z}_x$ from $S_{X_{z}}$ into $\mathcal{F}_{z}^2$ is uniformly continuous; and for $z'$ in the interior of $\mathbb{S}$, the map
$\phi_{z,z'}$ defined
by $\phi_{z,z'}(x)=F^z_x(z')$ is a uniform homeomorphism between the spheres of $X_z$ and $X_{z'}$.
\end{thm}

The following definitions and observations are used to clarify Theorem \ref{daher} and the work of Daher. Throughout this section
$(X_0, X_1)$ is a regular compatible couple of spaces.

\begin{fact}\label{compact} For any compact subset $K\subset \overset{\circ}{\mathbb{S}}$, there exists a constant $C > 0$ such that the norms
$\|.\|_{{\mathcal F}_z^2}$
and
$\|.\|_{{\mathcal F}_{z'}^2}$
are $C$-equivalent for all $z,z'$ in $K$.
\end{fact}

\pf It follows at once from the fact that the family of Poisson kernels $(\mu_{z})_{z \in K}$ consists of equivalent functions, in the sense that there are constants $c, C > 0$ such that $c \mu_{z} \leq \mu_{z'} \leq C \mu_{z}$ for all $z, z' \in K$. \pff

\begin{defi} The space $\mathcal{H}_{\infty}(X_0 + X_1)$ is the space of all functions $f : \mathbb{S} \rightarrow X_0 + X_1$ such that
\begin{enumerate}
    \item $f$ is analytic on the interior of $\mathbb{S}$,
    \item
        \[
        \|f\|_{\mathcal{H}_{\infty}} \coloneqq \sup\limits_{z \in \partial \mathbb{S}} \|f(z)\|_{X_0 + X_1} < \infty,
        \]
    \item $f$ satisfies the Poisson kernel representation
        \[
        f(z) = \int_{\partial \mathbb{S}} f(w) d\mu_{Re(z)}(w)
        \]
        for every $z$ in the interior of $\mathbb{S}$.
\end{enumerate}
\end{defi}

\begin{lemma}\label{43} Let $(X_0,X_1)$ be a regular couple of separable Banach spaces. Assume that $X_0^* \cap X_1^*$ is separable and that $(X_0^*,X_1^*)$ is regular. Let $f \in 
{\mathcal H}_\infty(X_0+X_1)$. Then $f$ admits non-tangential limits a.e. in $X_0+X_1$ equipped with the weak topology and coincides with this limit.\end{lemma}

\pf
Let $\varphi : \mathbb{D} \rightarrow \overset{\circ}{\mathbb{S}}$ be a conformal map, where $\mathbb{D}$ is the open disk. Then $f \circ \varphi$ is a function in the Hardy class $H^{\infty}(X_0 + X_1)$. Therefore, for any functional $\phi$ in $(X_0+X_1)^*=X_0^* \cap X_1^*$ the holomorphic function $\phi \circ f \circ \varphi$ admits non-tangential limits a.e.  By separability of $X_0^* \cap X_1^*$
we have non-tangential convergence a.e. of $\phi \circ f$ for $\phi$ in a countable dense family of functionals, and therefore non-tangential weak convergence a.e.
\pff

\

By convention when we write $f(1+it)$ for  $f \in {\mathcal H}_\infty(X_0+X_1)$ and $t \in \mathbb{R}$, it is understood that $t$ is such that the radial limit in $1+it$ exists and is equal to $f(1+it)$.

The previous lemma may be applied whenever we have a compatible couple $(X_0, X_1)$ of reflexive K\"othe function spaces with separable duals. Indeed, our assumptions imply that $X_0^*$ and $X_1^*$ are K\"othe function spaces as well. Let $f \in X_j^*$ be nonnegative and take a monotone sequence $(f_n)$ of simple functions such that $0 \leq f_n \leq f$ for every $n$ and $f_n \rightarrow f$ pointwise. Since $X_0^*$ and $X_1^*$ are reflexive Banach lattices they are $\sigma$-order continuous. This implies that $\|f - f_n\|_{X_j^*} \rightarrow 0$. Therefore, $X_0^* \cap X_1^*$ is dense in $X_0^*$ and in $X_1^*$.

Let us restate the definition of Daher's space $\mathcal{F}^{\infty} (X_0, X_1) $ but now calling attention to the space ${\mathcal H}_\infty (X_0 + X_1)$.

\begin{defi}[Daher]\label{def:F_infty} Let $(X_0,X_1)$ be a compatible couple of Banach spaces. Then $ {\mathcal F}^\infty(X_0,X_1)$ is the class of functions $f$ in ${\mathcal H}_\infty(X_0+X_1)$ such that
    \begin{enumerate}
        \item The functions $t \mapsto f(j + it)$, $j = 0, 1$, are measurable from $\mathbb{R}$ into $X_j$;
        \item \[
        \|f\|_{\mathcal{F}^{\infty}} := \max(\|f(i(\cdot))\|_{L^{\infty}(X_0)}, \|f(1 + i(\cdot))\|_{L^{\infty}(X_1)}) < \infty.
        \]
    \end{enumerate}
\end{defi}

\

Note that these functions are equal to their radial limits a.e. in the weak sense in $X_0+X_1$,  assuming $X_0^* \cap X_1^*$ is separable and the pair $(X_0^*,X_1^*)$ is regular.

\begin{prop}\label{weirdlimit} Let $(X_0,X_1)$ be a regular pair, and $0<\theta<1$.
Let $F \in {\mathcal F}^\infty$. Then 
$$\|F(1+it)\|_1=\lim_{z \rightarrow 1+it}\|F\|_{\mathcal F_z^2}$$
for $\mu_{\theta}^1$-almost every real $t$, including whenever $1+it$
is a point of continuity of
$z \mapsto \|F(z)\|_1$ on the right side of the border of $\mathbb S$. 
\end{prop}

\pf Recall the definition of the norm of ${\mathcal F}_{z}^2$:

$$\|F\|_{\mathcal{F}_{z}^2}^2 = \int\limits_{\partial \mathbb{S}} \|F(w)\|_z^2 d\mu_z(w) < \infty.$$ 

Note that $z \mapsto \|F\|_{\mathcal F_z^2}^2$ is actually the harmonic extension of the map $z \mapsto \|F(z)\|_z^2$ defined on the border.
In particular, we have radial limits a.e. for $\|F\|_{\mathcal F_z^2}$ when $z$ tends non tangentially to $z_0$ in the right vertical line, and the limit is then $\|F(z_0)\|_1$.
Under the continuity hypothesis, $\|F(1+it)\|_1$  is the radial limit in $1+it$ of this harmonic extension, hence the result follows. \pff

\

We shall observe next in Proposition \ref{unifunif} that if $X_0$ and $X_1$ are uniformly convex then we have uniform homeomorphisms, with modulus independent of $s \in (0,1)$, between the sphere of $X_s$ and the set of minimal functions equipped with the $\|.\|_{{\mathcal F}_s^2}$-norm.
Therefore limits considered in the previous and in the present sections, such as  $\lim_{s \rightarrow 1} \|F(s)\|_s$ and $\lim_{s \rightarrow 1}\|F\|_{\mathcal F_s^2}$, are actually related, when $F$ is a difference of minimal functions.

\subsection{Uniform continuities}

Recall Daher's result (Theorem \ref{daher}): if $(X_0,X_1)$ is a regular scale of reflexive spaces and $X_0$ or $X_1$ is strictly convex, then the map
$$\Gamma_z: S_{X_z} \rightarrow {\mathcal F}_z^2$$
is well defined. Moreover, if $X_0$ or $X_1$ is uniformly convex, then $\Gamma_z$ is uniformly continuous for every $z$ in the interior of the strip $\mathbb S$. 

Let us note the following example which shows that uniform convexity is indeed required for the result of the previous paragraph.

\begin{example} 
The scale $(\ell_1,\ell_\infty)$ shows that the above does not hold if one removes the uniform convexity assumption. Indeed, in this case note that $X_\theta=\ell_{\frac{1}{1-\theta}}$, and we have the formula
$$\Gamma_\theta(x)(z)=x^{\frac{1-z}{1-\theta}}$$
if $x \geq 0$ of norm $1$ in $X_\theta$,
and $$\Gamma_\theta(x)=sgn(x) |x|^{\frac{1-z}{1-\theta}}$$
for $x \in S_{X_\theta}$. Let $0 < \alpha < 1$, $1 / 2 < \theta < 1$ and $x = (x_n), y = (y_n) \in S_{X_\theta}$ such that $x > 0$, $x_1 < \alpha$, $y_1 = -x_1 $, and $y_n = x_n$ for all $n\geq 2$. Then 
$$\|x-y\|_{\theta}=\| (sgn(x)-sgn(y))|x| \|_{\ell_{\frac{1}{1-\theta}}} \leq 2 \alpha,$$
however 
\begin{align*}
\|\Gamma_\theta(x)-\Gamma_\theta(y)\|_{\mathcal{F}_\theta^\infty} &= \int_{t \in \R} \|\Gamma_\theta(x)(1+it)-\Gamma_\theta(y)(1+it)\|_\infty d\mu_\theta^1(t) \\ 
&=\int_{t \in \R} \|sgn(x)|x|^{\frac{it}{1-\theta}}-sgn(y)|y|^{\frac{it}{1-\theta}}\|_\infty d\mu_\theta^1(t) \\
&= \int_{t \in \R} \| sgn(x)-sgn(y)\|_\infty d\mu_\theta^1(t)=\theta  \| sgn(x)-sgn(y)\|_\infty > 1,
\end{align*}
proving the claim. 
\end{example}

\

In addition to the uniform continuity of $\Gamma_z: S_{X_z} \rightarrow {\mathcal F}_z^2$, whenever $z \in \overset{\circ}{\mathbb{S}}$, we obtain the following strengthening of the result: 

\begin{prop}\label{unifunif}
Let $(X_0,X_1)$ be a regular couple of interpolation with $X_0$ and $X_1$ uniformly convex. Then the modulus of continuity of $\Gamma_{z}$ is uniform for $z$ in $\overset{\circ}{\mathbb{S}}$.
\end{prop}
\pf
First note that by vertical symmetry, it is enough to prove uniformity of $\Gamma_\theta$ for $0<\theta<1$.
Indeed we have the relation
$$
\Gamma_{\theta+it}(x)(z)=\Gamma_{\theta}(x)(z-it)
$$
and therefore, if $$v_t: {\mathcal F}_\theta^2 \rightarrow {\mathcal F}_{\theta+it}^2$$ denotes the isometric map 
defined by $v_t(F)(z)=F(z-it)$, we have 
$\Gamma_{\theta+it}=v_t \circ \Gamma_{\theta}$. So the moduli of continuity of $\Gamma_\theta$ and $\Gamma_{\theta+it}$ are equal.

Fix $j = 0, 1$, and note that the spaces $L_2(\mu_{\theta}^{j})$ have the same density character for every $\theta$, so they are all isometric. By tensoring these isometries with the identity map, we get that the spaces $L_2(\mu_{\theta}^{j}, X_j)$ are all isometric as well, so they have the same modulus of convexity.

Daher shows that the modulus of continuity of $\Gamma_{\theta}$ may be bounded in terms of the modulus of convexity $\delta_0$ of $L_2(\mu_{\theta}^{0}, X_0)$ and $\delta_1$ of $L_2(\mu_{\theta}^{1}, X_1)$.
Indeed, if $\delta = \min\{\delta_0, \delta_1\}$, Daher shows that
\[
\delta \left (\frac{1}{2}\|\Gamma_{\theta}(x) - \Gamma_{\theta}(y)\|_{\mathcal{F}_{\theta}^2} \right ) \leq \|x - y\|_{\theta}
\]

\noindent for every $\theta \in (0, 1)$ and $x, y \in S_{X_{\theta}}$. Let $\beta(t) = \sup\{u \geq 0 : \delta(u) \leq t\}$. It follows that
\[
\frac{1}{2} \|\Gamma_{\theta}(x) - \Gamma_{\theta}(y)\|_{\mathcal{F}_{\theta}^2} \leq \beta(\|x - y\|_{\theta})
\]

Now, since $\delta$ is a continuous non decreasing function such that $\lim_{t \rightarrow 0} \delta(t) = 0$, it follows that $\lim_{t \rightarrow 0} \beta(t) = 0$, and therefore $\Gamma_{\theta}$ is uniformly continuous.

Note, however, that the previous estimates do not depend on $\theta$. The result follows.
\pff

\

Note that the evaluation map $F \mapsto F(z)$ is by definition bounded of norm $1$ from ${\mathcal F}^2_z$ onto $X_z$.
Nevertheless, there is no reason to assume that the evaluation maps $F \mapsto F(\theta)$ are uniformly bounded when $\theta$ tends to the border and the norm of $F$ is with respect to a fixed ${\mathcal F}_{x_0}^2$.

\

From Proposition \ref{unifunif} we deduce the following technical result, that will be of use later on.  Let $\beta$ be the map satisfying $\lim_{t \rightarrow 0}\beta(t)=0$ defined in the proof of Proposition \ref{unifunif}.

\begin{lemma}\label{new}
Let $(X_0,X_1)$ be a regular compatible couple of uniformly convex spaces. Fix $0<\theta<1$. Let $\varepsilon>0$. If $x,y \in S_{X_\theta}$
then there exists $T \subseteq \mathbb{R}$ with full $\mu_\theta^1$-measure such that
for any $t \in T$,
$$\limsup_{s \rightarrow 1^{-}} \|\Gamma_{\theta}(x)(s+it)-
\Gamma_{\theta}(y)(s+it)\|_s$$ 
$$\leq \|\Gamma_{\theta}(x)(1+it)-\Gamma_{\theta}(y)(1+it)
\|_1 $$ 
$$\leq 
2\beta(\limsup_{s \rightarrow 1^{-}}  \|\Gamma_{\theta}(x)(s+it)-
\Gamma_{\theta}(y)(s+it)\|_s).$$
\end{lemma}

\pf It follows from Proposition \ref{weirdlimit} that there is a measurable set $T \subseteq \mathbb R$ with full $\mu_\theta^1$-measure such that for all $t \in T$ $$\|\Gamma_{\theta}(x)(1+it)-\Gamma_{\theta}(y)(1+it)
\|_1 =
\lim_{s \rightarrow 1^{-}}
\|\Gamma_{\theta}(x)-
\Gamma_{\theta}(y)\|_{{\mathcal F}_{s+it}^2}.$$

The left inequality then follows from the fact that
$$\|(\Gamma_{\theta}(x)-
\Gamma_{\theta}(y))(s+it)\|_s=
\|(\Gamma_{\theta}(x)-
\Gamma_{\theta}(y))(s+it)\|_{s+it} \leq
\|\Gamma_{\theta}(x)-
\Gamma_{\theta}(y)\|_{{\mathcal F}_{s+it}^2}.
$$
The right inequality follows from the fact that
$$\|\Gamma_{\theta}(x)-
\Gamma_{\theta}(y)\|_{{\mathcal F}_{s+it}^2}
\leq 2\beta(\|(\Gamma_{\theta}(x)-
\Gamma_{\theta}(y))(s+it)\|_{s+it});$$
this fact is a consequence of the estimate for the modulus of $\Gamma_s$ in
the proof of Proposition
 \ref{unifunif}, plus the fact that $\Gamma_s$ and $\Gamma_{s+it}$ have the same modulus). Namely
 $$\|\Gamma_{s+it}(x')-
\Gamma_{s+it}(y')\|_{{\mathcal F}_{s+it}^2}
\leq 2\beta(\|x'-y'\|_{s+it}),$$
 applied to  $x'=\Gamma_{\theta}(x)(s+it)$ and
 $y'=\Gamma_{\theta}(y)(s+it)$, noting that 
 $\Gamma_{s+it}(x')=\Gamma_{\theta}(x)$
 and $\Gamma_{s+it}(y')=\Gamma_{\theta}(y)$.
  \pff

\

Note that the result of Lemma \ref{new} easily follows from Corollary \ref{corocoro} in the case of an interpolation scale of spaces with a common $1$-Schauder basis, or more generally a common 1-FDD. However, we prefer to follow as much as possible in the vein of Daher's results, which do not require any basis structure.

\subsection{Vertical maps}

A second ingredient to obtain properties of interpolation scales close to the border will be the definition and study of vertical maps.
In this section 
 $(X_0,X_1)$ is a regular couple, one of them is uniformly convex, and $0<\theta<1$ is fixed.
Recall that we denote as $\phi_{z,z'}$ the uniform homeomorphism between the spheres of $X_z$ and $X_{z'}$ defined by
$$\phi_{z,z'}(x):=(\Gamma_z(x))(z').$$

Note in particular that $\phi_{\theta+it,\theta+it'}$ is a uniform homeomorphism on $S_{X_\theta}$. By vertical symmetry we have the identities
$$\Gamma_{\theta}(\Gamma_{\theta}(x)(\theta + is))(z) = \Gamma_{\theta}(x)(z + is)$$ and
$$\phi_{\theta+it,\theta+it'}=\phi_{\theta,\theta+i(t'-t)}.$$
Therefore we denote
\begin{defi}
Let $\{\phi_\theta^t : t \in \R\}$  be the  $1$-parameter group of uniform homeomorphisms on $S_{X_\theta}$ defined by
 $$\phi_\theta^t:=\phi_{\theta,\theta+it} : S_{X_\theta} \rightarrow S_{X_\theta}.$$   We shall call such uniform homeomorphisms the {\em vertical homeomorphisms} on $S_{X_\theta}$.
\end{defi}

\

Vertical maps, in spite of being defined on a fixed $X_\theta$, give a lot of information on the whole scale through the map $\Omega_\theta$ discussed in the introduction, using as an advantage that they are defined as non-linear operators on the same space, as opposed to maps $\phi_{\theta,\theta'}$ which act from (the sphere of) $X_\theta$ to $X_{\theta'}$. For example, in \cite{Stability}, it is proved that, under reasonable hypotheses, if the vertical maps on $S_{X_\theta}$ are constant for some fixed $\theta$ (which corresponds to $\Omega_\theta=0$), then the scale is trivial, i.e. $X_0=X_1$ isometrically. Some results are also obtained in \cite{Stability} when $\Omega_\theta$ is linear, also see \cite{Suarez} for the case when $\Omega_\theta$ is bounded.

The vertical maps can be explicitly computed in very specific cases. In the case of the $\ell_p$-scale they can be seen as a vertical version of the Mazur maps, as described below.

\begin{example}\label{verticalforlp} Consider $X_0=\ell_\infty$ and $X_1=\ell_q, 1 \leq q <+\infty$. Then given $\theta \in (0,1)$, we have
$$\phi_\theta^t(x)=x |x|^{it/\theta}.$$
Indeed, note that ${X_\theta} = \ell_p$ where  $p=q/\theta$. Thus, given $x \in S_{X_\theta} $, note that $F_x^\theta(z)=sgn(x) |x|^{ z /\theta}$, and the vertical maps are 
$$\phi_\theta^t(x)=F_x^\theta(\theta+it)=sgn(x) |x|^{1+it/\theta }=x |x|^{it/\theta}.$$

\end{example}

\

It is interesting to have a look at the moduli of uniform continuity of those maps.

Take $f(x) = x \left|x\right|^{\frac{it}{\theta}}$. Then
we can write
$$f(y)-f(x)=(y-x)|x|^{\frac{it}{\theta}}+y(|y|^{\frac{it}{\theta}}-|x|^{\frac{it}{\theta}})$$
The maps $r \rightarrow r^{it/\theta}$ are uniformly continuous on $[0,1]$ with a modulus of uniform continuity $\rho$ for all $t/\theta \leq M$, given $M>0$. 
We deduce that whenever $\|x-y\|_\theta \leq \alpha$,
$$\|f(y)-f(x)\|_\theta \leq \|y-x\|_\theta+\rho(\alpha).$$

In particular the vertical maps have a common modulus of uniform continuity whenever $t/\theta$ belongs to a bounded set.

\

When $\theta$ tends to $1$, note that the above map seems to tend naturally to
the map $x \mapsto x|x|^{it}$, which is again a uniform homeomorphism on $S_{X_1}=S_{\ell_q}$. In what follows we shall investigate conditions for this situation to happen, which in turn gives favorable conditions to push uniform homeomorphism from the interior of the scale to the border.

\subsection{The limit map}

By Lemma \ref{limit},
if, for example, a regular scale of spaces with a common $1$-Schauder basis is continuous to the right, then we may define a natural map $\phi_{\theta,1}$ from $S_{X_\theta} \cap c_{00}$ to $S_{X_1} \cap c_{00}$, by 
$$\phi_{\theta,1}(x)=\lim_{s \rightarrow 1} \phi_{\theta,s}(x).$$  In the best case scenario, this map could provide 
\begin{itemize}
 \item 
a uniform homeomorphism between the spheres of $X_\theta$ and $X_1$, or,
\item in the case of interpolation of $X_1$ with its antidual, a map witnessing that $X_1$ has property (H).
\end{itemize}

Even under weaker continuity hypothesis, a limit map may be defined and studied. In what follows we assume a regular pair $(X_0,X_1)$ of uniformly convex spaces which has separable intersection space, as in Lemma \ref{43}. Note that the dual scale is regular as well, and $X_0$ and $X_1$ are separable. We shall call such compatible couples \emph{adequate}.

\begin{defi}  \label{def:mazurmap} Let $(X_0,X_1)$ be an adequate pair.
We let $\Delta_\theta \subset S_{X_\theta}$ be the set of $x$ such that $F_x^\theta(s)$ tends weakly to $F_x^\theta(1)$ when $s \rightarrow 1$ (in $X_0+X_1$) and such that $F_x^\theta(1)$ belongs to $S_{X_1}$.
We define the \emph{limit Mazur map} $\phi_{\theta,1}: \Delta_\theta \rightarrow X_1$ by
$$\phi_{\theta,1}(x):=F_x^\theta(1)=w-\lim_{s \rightarrow 1-} \phi_{\theta,s}(x),$$ 
and we denote $$S_1:=\phi_{\theta,1}(\Delta_\theta)$$ for some, and therefore for any $0<\theta<1$.
\end{defi}

We also denote $S_1^{\prime}$ to be the set of elements of $X_0+X_1$ of the form $\phi_{\theta,1}:=w-\lim_s F_x^\theta(s)$ for some $\theta$ and $x \in S_{X_\theta}$. Therefore
$S_1=S_1^{\prime} \cap S_{X_1}$.

\begin{lemma}\label{nonconstant} Let $(X_0,X_1)$ be a regular couple. A non-constant holomorphic function in $H_\infty(X_0+X_1)$ cannot be constant on a set of positive measure on the boundary. In particular, it assumes uncountably many values on the right boundary.
\end{lemma}

\pf Otherwise the composition of such $F$ with any dual functional will be constant on a set of positive measure on the border. Since by \cite{D} (Theorem 2.2 and the comment at the end of the page), 
for $f \in H_\infty$, $f(z)$ is identically $0$ if it vanishes on a set of positive measure on the boundary, then we will have that the composition of $F$ with any dual functional is actually constant, whereby $F$ is itself constant.  
\pff

\begin{lemma} The set $S_1^{\prime}$ contains all radial limits at all horizontal levels of all minimal functions. For each minimal function these radial limits are a.e. in $S_{X_1}$. The set $S_1$ is uncountable. \end{lemma}

\pf Fix $\theta$. Note that given $x \in S_{X_\theta}$ and $t$ such $F_x^\theta$ has radial limit $y$ in $1+it$, this means that $F_{x'}^{\theta+it}$ has radial limit in $1+it$, where $x'=F_{\theta,\theta+it}(x)$. Then by vertical symmetry this means that $F_{x'}^{\theta}$ tends to $y$ in $1$. The comment after Definition \ref{def:F_infty} and Daher's Theorem \ref{thm:Daher} ensure that almost every time $y \in S_{X_1}$.

If $X_0=X_1$ then the scale is trivial and it is clear that $S_1=S_{X_1}$, so the result trivially holds. If they are different then there must exist some non-constant minimal function and we are done by the first part of the proof. \pff

\

\begin{defi} Let $(X_0,X_1)$ be a regular couple of uniformly convex spaces.
For $x \in S_{X_\theta}$, let $O_\theta(x) \subset S_{X_\theta}$ be the orbit of $x$ under the action of the vertical homeomorphisms. In other words, 
$$O_\theta(x):=F_x^{\theta}(\theta+i\R)=\{\phi_{\theta}^t(x) : \, t \in \R\}.$$
\end{defi}

We say that a property $P$ of elements of $S_{X_\theta}$ holds a.e. (resp. almost nowhere) on an orbit $O_{\theta}(x)$ if $\{t: \phi_\theta^t(x) {\rm \  satisfies\ } P \}$ has full measure  (resp. null measure) with respect to the measure $\mu_\theta^1$ on $\R$.

\begin{lemma}\label{lem:Mazur_orbit} Let $(X_0,X_1)$ be an adequate compatible couple. Let $x \in S_{X_\theta}$. Then
\begin{itemize}
\item[(a)] $\phi_{\theta,1}(y)$ is defined and with value in $S_{X_1}$ a.e.\ on the orbit $O_{\theta}(x)$,
 
\item[(b)] if
$F_x^\theta$ is not constant and  $x' \in S_1$, then $\phi_{\theta,1}(y)=x'$ holds almost nowhere on $O_{\theta}(x)$. 
\end{itemize}
\end{lemma}

\pf \emph{(a)} This a consequence of  Lemma \ref{43} applied to the map $F_x^\theta$.
\medskip

\emph{(b)} If the measure of $\{t: \phi_{\theta}^t(x) \in \phi_{\theta,1}^{-1}(y)\}$ was not 0, it would follow that $F_x^{\theta}$ is a constant function, by Lemma \ref{nonconstant}.

\pff

\begin{fact}\label{lem:continuity}
Let $F \in \mathcal{F}_{\theta}^{1}$, where $\theta \in (0, 1)$. Then $t \mapsto F(\theta + it)$ is a continuous function from $\mathbb{R}$ into $X_{\theta}$.
\end{fact}
\pf
For simplicity, we prove continuity at $0$. Let $j = 0, 1$. Since $\nu_{\theta, j}$ is an atomless regular measure, for every simple function $G:\mathbb R\rightarrow X_j$ of the form $G=\sum_{k=1}^m \chi_{E_k} x_k$ with $E_k\subset\mathbb R$ pairwise disjoint compact sets and $x_k\in X_j$, and for $|t|<\min\{d(E_k,E_l):k\neq l\}$, we have
\begin{equation}\label{eq:translation on simple goes to 0}
\int_{\mathbb R}\|G(s)-G(s+t)\| d\nu_{\theta, j}(s)=\sum_{k=1}^m\nu_{\theta, j}\Big(E_k\Delta (E_k-t)\Big)\|x_k\|\underset{t\rightarrow 0}\longrightarrow 0,
\end{equation}
where $E_k\Delta (E_k-t)$ denotes the symmetric difference of the set $E_k$ and $E_k-t=\{s\in \mathbb R: s+t\in E_k\}$. 

Note that since $\mu_{\theta, j}$ is regular, the functions of the form $G$ are dense in the Bochner space $L_1(\mu_{\theta}^{j},X_j)$. Hence, given $F\in \mathcal{F}_{\theta}^{1}$, taking simple functions of the above form which approximate $F(j+i\cdot)$ we get that
$$
\int_{\mathbb R}\|F(j+is)-F(j+i(s+t))\| d\mu_{\theta}^{ j}(s)\underset{t\rightarrow 0}\longrightarrow 0.
$$

Now the result follows from the definition of the norm in $X_{\theta}$.
\pff

\begin{coro}\label{denseness} 
The set $\Delta_\theta$ is dense in $S_{X_\theta}$ \end{coro}

\pf For all $x \in S_{X_\theta}$, $\phi_{\theta,1}$ is defined a.e.\ on $O_{\theta}(x)$ with values  in $S_{X_1}$; 
Therefore there are points of $\Delta_\theta$ arbitrarily close to $x$.\pff

\

Let $X$ be a Banach space of finite dimension $n$. Let us consider the map that sends $x \in S_{\ell_2^n}$ into $\frac{x}{\|x\|_X} \in S_X$. This map allows us to transport the canonical surface measure on $S_{\ell_2^n}$ to $S_X$. We will always assume that if $X$ is a finite dimensional Banach space then $S_X$ is equipped with that measure.

\begin{prop}\label{props:fin_dim_measure_1} Assume $X_0, X_1$ are finite-dimensional. Then $\Delta_\theta$ is a subset of measure $1$ of $S_{X_\theta}$.
\end{prop}

\pf We know that 
$$\forall x \in S_{X_\theta}, \forall^* t \in \R, \phi_{\theta}^t(x) \in \Delta_\theta.$$
Therefore by Fubini's theorem
$$\forall^* t \in \R, \forall^* x \in S_{X_\theta}, \phi_{\theta}^t(x) \in \Delta_\theta.$$
By vertical symmetry and since $\phi_{\theta}^t$ is a homeomorphism for fixed $t$, we deduce the result from the case $t=0$.
\pff

\medskip
Summing up we have the following:

\begin{prop}\label{props:Delta_theta_in_S_theta} Let $(X_0,X_1)$ be a regular couple of uniformly convex spaces which are finite dimensional or with a common $1$-Schauder basis $(e_n)$.  Then: 
\begin{itemize} 
\item[(1)] $\Delta_\theta \cap [e_1,\ldots,e_n]$ has measure $1$
in $S_{X_\theta^n}$ ($0 < \theta < 1$).
\item[(2)] $\Delta_\theta \cap c_{00}$ is equal to $S_{X_\theta} \cap c_{00}$ if the scale is continuous to the right, i.e., we have a limit map $\gamma_{\theta,1}$ from $S_{X_\theta} \cap c_{00}$ into $S_{X_1} \cap c_{00}$.
\end{itemize}
\end{prop}
\pf
(1) The finite dimensional case was done in Proposition \ref{props:fin_dim_measure_1}. For the $1$-Schauder basis case, we have a common norm $1$ projection from $X_j$ onto $X_j^n$ ($j = 0, 1$). Therefore $X_{\theta}^n = (X_0^n, X_1^n)$ isometrically, if we denote the span of $e_1, \cdots, e_n$ in $X_{\theta}$ by $X_{\theta}^n$, and we are back to the finite dimensional case.

(2) Let us consider first the finite dimensional case and take $x \in S_{X_\theta}$. Since the scale is continuous to the right, the minimal function $F_x^{\theta}$ is continuous on $\{z \in \mathbb{S} : 0 < Re(z) \leq 1\}$, and therefore there exists $\phi_{\theta, 1}(x)$. We only have to check that it belongs to $S_{X_1}$. Since $\|.\|_{X_1}$ is equivalent to $\|\cdot\|_{X_0 + X_1}$ the map $t \mapsto \|F_x^{\theta}(1 + it)\|_1$ is continuous. By Lemma \ref{lem:Mazur_orbit}(a) we must have that $\phi_{\theta, 1}(x) \in S_{X_1}$. Therefore $S_{X_{\theta}} = \Delta_{X_{\theta}}$. Now the $1$-Schauder basis case follows directly.

\pff

In case (2) we can ask when this map is an homeomorphism, a witness for property (H), or a uniform homeomorphism. 
Before this, we discuss the possible behaviour of the vertical orbits.

\

\subsection{On the vertical orbits}

In general, the behaviour of the limit map seems to depend on the vertical orbits. We therefore give an analysis of those based on a classification of points in $X_\theta$.

\begin{defi} We say that a point $x$ of $S_{X_\theta}$ is $T$-periodic for some $T>0$ (resp. periodic, aperiodic, singular), if the map
$t \mapsto \phi_\theta^t(x)$ is $T$-periodic (resp. periodic, aperiodic, constant). 
\end{defi}

The previous definition is related to the value of $\inf \{t>0: F_x^\theta(t)=x\}$, which can be either $0$, positive or infinite. This corresponds to points of $S_{X_\theta}$ which are either singular, periodic non singular, or aperiodic. Each of the orbits is respectively either a point, homeomorphic to the unit circle, or homeomorphic to $\R$.

\

Some observations are in order: first, the sets of singular points, $T$-periodic points and aperiodic points are stable under multiplication by $\lambda$, with $\left|\lambda\right| = 1$. Also, 

\begin{prop} Assume we have a regular scale of uniformly convex spaces. Then the sets of singular points and of $T$-periodic points are closed. \end{prop}

\pf A point $x \in S_{X_\theta}$ is singular if and only if $\|x\|_0=\|x\|_1=1$.
Since the map taking $x$ to its minimal function in $\theta$ is uniformly continuous in the appropriate norm, 
if $x_n$ is a sequence of singular points converging to $x$, then $F_{x_n}^\theta$ converges to $F_x^\theta$. Since all functions of the sequence are constant, the limit is constant as well, so $x$ is singular. The same argument holds for $T$-periodic.
\pff

\begin{example}
Consider the compatible couple $(\ell_p , \ell_q )$ for $1 < p , q < \infty$, and let $1 \slash p_{\theta} = (1 - \theta ) \slash p + \theta \slash q$. For $x = (x_n) \in S_{\ell_{p_{\theta}}}$, consider $\psi_x : \R \longrightarrow S_{\ell_{p_{\theta}}}$ given by $\psi_x (t) \coloneqq \phi_{\theta}^t (x) = (x_n | x_n|^{ita} )$, where $a = p_{\theta} (1 \slash p - 1 \slash q)$. Note that $\psi_x$ is constant if and only if for each $n \in \N$ we have that $x_n = 0 $ or $|x_n| =1 $. 

Let $A = \{ n \in \N :  x_n \neq 0 , |x_n| \neq 1 \}$. Let us see that $\psi_x$ is $T$-periodic if and only if there is $C > 0$ and for each $n \in A$ there is $m(n) \in \Z$ such that $m(n) \slash \ln (|x_n| ) = C$.  In this case, the period $T$ of $\psi_x$ is given by $T = 2 \pi C \slash a$. 
Indeed, if $\psi_x$ is $T$-periodic, we have that for every $t \in \R$ 
$$ x_n | x_n|^{ita} = x_n | x_n|^{ia(t+T)}, \ \forall n \in \N  .$$
In particular, $x_n | x_n|^{ita} = x_n | x_n|^{ia(t+T)}, \ \forall n \in A  $. Thus, $|x_n|^{iaT} =1 $ for all $n \in A$, which implies that for each $n \in A$ there is $m(n) \in \Z$ such that $T = 2 \pi m(n) \slash a \ln (|x_n|) $.
\end{example}

\

Note that if the set of singular points is dense in $S_{X_{\theta}}$, then $X_0 = X_1$. 
However, this cannot be much weakened as the following shows:

\begin{example}
This example shows that even if the set of singular points is dense in some open set of the sphere of $X_{\theta}$, it does not follow that $X_0 = X_1$. We give a geometric argument of why such an example must exist, but it is important to remember that our Banach spaces are complex. Consider $X_0 = \ell_2^2$ and the two dimensional Banach space $X_1$ which unit ball is similar to:

\begin{figure}[H]
\includegraphics[width=8cm]{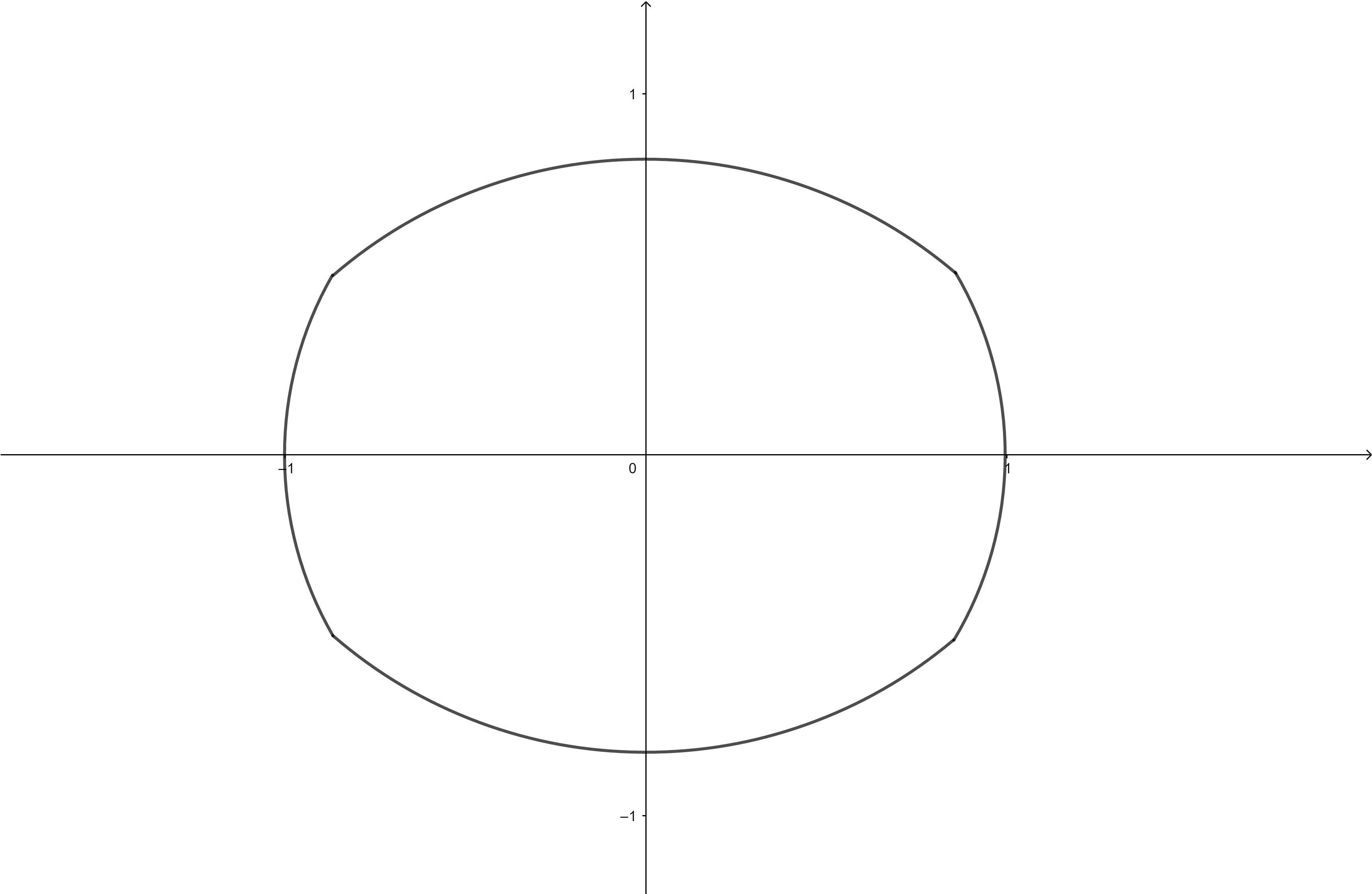}
\centering
\end{figure}

The left and right sides above agree with the sphere of $\ell_2^2$. So $S_{X_0} \cap S_{X_1}$ is a nonempty open subset of $X_{\theta}$ for every $0 < \theta < 1$. Also, any $x \in S_{X_0} \cap S_{X_1}$ is a singular point.
\end{example}

Finally, we observe that the existence of periodic points imposes heavy restrictions on the minimal functions.

\begin{prop}\label{cndn} Assume $x \in S_{X_\theta}$ is $T$-periodic and let $\phi$ be its norming functional. Then
\begin{itemize}
    \item[(1)] $\phi$ is $T$-periodic as well;
    \item[(2)]
    if $(X_0,X_1)$ are finite dimensional, then we have the formulas
$$F_x^\theta(z)=\Sigma_n c_{n} {\rm exp}( 2\pi n (z-\theta)/T ),$$
$$F_\phi^\theta(z)=\Sigma_n d_{n} {\rm exp}( 2\pi n (z-\theta)/T ),$$
with $\sum_n c_n$ and $\sum_n d_n$ normally convergent in $X_\theta$ and
$X_\theta^*$ respectively,
$1=\sum_n \langle d_{-n},c_n\rangle$
and for all $k \neq 0$,
$0=\sum_{n+m=k} \langle d_{n},c_{m}\rangle.$
\end{itemize}

\end{prop}

\pf (1) Regarding $F_\phi^\theta$, since
$$1=\langle\phi,x\rangle=\int_{\partial \mathbb{S}} \langle F_\phi^\theta(z),F_x^\theta(z)\rangle d\nu_{\theta} (z),$$
necessarily $F_\phi^\theta(z)$ is the unique norming functional of $F_x^\theta(z)$, so whenever $F_x^\theta(z)=x$ then 
$F_\phi^\theta(z)=\phi$. This proves that $\phi$ is $T$-periodic whenever $x$ is.

(2)
By the characterization of periodic functions on page 182 of Gamelin \cite{Ga}, it follows that we have the formula above for $F_x^\theta$ and for $F_\phi^\theta$, associated to normally convergent series $\sum_n c_n$ and $\sum_n d_n$.

From this it follows that
$$1=\langle F_\phi^\theta(z), F_x^\theta(z)\rangle= \Sigma_{m,n} \langle d_m, c_n \rangle {\rm exp}((m+n)(2\pi (z-\theta)/T)).$$
Which means (taking $z = \theta$)
$$1=\sum_n \langle d_{-n},c_n\rangle,$$
and for all $k \neq 0$,
$$0=\sum_{n+m=k} \langle d_{n},c_{m}\rangle.$$

\pff

The example of scales $\ell_p$-spaces show that periodic points do appear in certain special cases.  The following is an easy exercise, and indicates that lattice structure seems to be a condition to give periodic points; it is unclear however if this may happen concretely
apart from the case of the scale of $\ell_p$-spaces.

\begin{example} Assume $c_n$ is a common $1$-unconditional basis for an interpolation scale with $\sum_n \|c_n\|_\theta<+\infty$, and such that $x:=\sum_n c_n$ belongs to $S_{X_\theta}$. Assume furthermore that $$\|\Sigma_n c_n {\rm exp}(n(1-\theta))\|_1=1$$
and
$$\|\Sigma_n c_n {\rm exp}(-n\theta)\|_0=1.$$
Then $x$ is $2\pi$-periodic, and 
$F_x^\theta(z)=\Sigma_n c_n {\rm exp}(n(z-\theta))$.
\end{example}

\medskip

\section{Uniform scales}\label{sec:uniform_scales}

In this section we consider an ultimate list of ``internal" properties of interpolation scales trying to mimic extrapolation, but which are formally weaker. We shall call them ``uniformities", and they will rely heavily on the use of the vertical homeomorphisms considered in the previous section.

Before giving the definition, we consider general technical lemmas relating values of analytic functions on the border and in the interior of the strip.

\subsection{Interior versus border of the strip}

For notational ease and to avoid unrelevant technicalities, we let $\gamma_\theta^1$ denote the normalization of the measure $\mu_\theta^1$ in $\R$, i.e.
$$\gamma_\theta^1(A):=\mu_\theta^1(A)/\theta,$$
so that $\gamma_\theta^1$ is a probability measure on $\R$.
Recall  that, as already used in Fact \ref{compact}, for fixed $\theta$, the measures $\gamma_\theta^1$, $\mu_\theta^1$, $\gamma_{1/2}^1$ and $\mu_{1/2}^1$ are all equivalent, so that the notions of a.e. and of sets of positive measure coincide for all of them, as well as sets of small measure up to a multiplicative constant. The measure $\gamma_\theta^1$ is simpler to use in the initial part of this section, and then we shall switch to $\gamma_{1/2}^1$, which makes more sense for global properties of a scale, in the next section.

\begin{lemma}\label{classicalinterpolationestimate}
Let $(X_0,X_1)$ be a regular compatible couple.
Let 
$0<\theta<1$ and $F \in {\mathcal F}_\theta^1$.
Then $$\|F(\theta)\|_\theta \leq  \Big(\int_{\R} \|F(it)\|_0 d\gamma_{\theta}^0(t)\Big)^{1-\theta} \Big(\int_{\R} \|F(1+it)\|_1 d\gamma_{\theta}^1(t)\Big)^\theta.$$
\end{lemma}

\pf By definition,
$$\|F(\theta)\|_\theta \leq 
(1-\theta) \int_{\R} \|G(it)\|_0 d\gamma_{\theta}^0(t)
+\theta \int_{\R} \|G(1+it)\|_1 d\gamma_{\theta}^1(t)
$$
whenever $G\in {\mathcal F}_\theta^1$ satisfies $G(\theta)=F(\theta)$. Using $G$ defined by the formula
$G(z)=e^{K(z-\theta)} F(z)$, with $K \in \R$, we deduce
$$\|F(\theta)\|_\theta \leq 
(1-\theta)e^{-K\theta}(\int_{\R} \|F(it)\|_0 d\gamma_{\theta}^0(t))
+\theta e^{K(1-\theta)} (\int_{\R} \|F(1+it)\|_1 d\gamma_{\theta}^1(t)).
$$
Minimizing over $K$ gives the expected result.
\pff

\begin{lemma}\label{border}  Let $(X_0,X_1)$ be a regular compatible couple.
Let 
$0<\theta<1$ and let $\varepsilon>0, c>0$.
\begin{itemize}
    \item[(1)] There exists $\alpha>0$ such that if $\|F\|_{{\mathcal F}_\theta^2}<\alpha$, then $$\gamma_{1/2}^1(\{ t \in \R: \|F(1+it)\|_1 \leq \varepsilon\} )\geq 1-\varepsilon.$$
    \item[(2)] There exists $\alpha>0$ such that whenever
    $\|F\|_{{\mathcal F}_\infty}\leq c$ and
     $\gamma_\theta^1(\{t \in \R: \|F(1+it)\|_1 \leq \alpha\}) \geq 1-\alpha/c$, then $\|F(\theta)\|_{\theta}<\varepsilon$.
\end{itemize}
\end{lemma}

\pf (1) 
This is a consequence of the integral definition of the norm of $F$. We have:
\begin{eqnarray*}
\|F\|_{\mathcal{F}_{\theta}^2}^2 & = & \int_{\partial \mathbb{S}} \|F(z)\|_z^2 d\mu_{\theta}(z) \\
    & = & \int_{Re(z) = 0} \|F(z)\|_0^2 d\mu_{\theta}(z) + \int_{A} \|F(z)\|_z^2 d\mu_{\theta}(z) + \int_{\mathbb{S}_1 \setminus A} \|F(z)\|_z^2 d\mu_{\theta}(z)
\end{eqnarray*}
for every measurable set $A \subset \mathbb{S}_1 \coloneqq \{z  \in \C: Re(z) = 1\}$. Take $A$ to be the set where $\|F(z)\|_1> \varepsilon$. If $\|F\|_{\mathcal{F}_\theta^2} < \alpha$ then
\begin{eqnarray*}
\mu_{\theta}(A) \epsilon^2 \leq \|F\|_{\mathcal{F}_{\theta}^2}^2 < \alpha^2
\end{eqnarray*}
which implies that the $\gamma_{\theta}^1$-measure of the complement of $A$ in $\mathbb{S}_1$ is at least $1 - \frac{\alpha^2}{\theta \varepsilon^2}$. So it is enough to take $\alpha = \varepsilon^{3/2} \theta^{1/2}$.

(2) Let  $A \subset \R$ be the set where $\|F(1+it)\|_1 \leq \alpha$. We have:
$$\int_{\R} \|F(1+it)\|_1 d\gamma_{\theta}^1(t) \leq \alpha \gamma_{\theta}^1(A) +c(1-\gamma_{\theta}^1(A)) \leq \alpha +c(1-\gamma_{\theta}^1(A)) \leq 2\alpha.$$
The classical interpolation estimate from Lemma \ref{classicalinterpolationestimate} may then be used to obtain that
 $\|F(\theta)\|_\theta \leq c^{1-\theta} (2\alpha)^\theta$ which is at most $\varepsilon$ if $\alpha$ was well-chosen.
\pff

 \begin{lemma}\label{55}
 Let $(X_0,X_1)$ be an optimal compatible couple. Then
 \begin{itemize}
 \item[(1)]
 For every $\varepsilon>0$, there exists $\alpha>0$ such that whenever $x,y \in S_{X_\theta}$ satisfy
 $$\gamma_\theta^1(\{t \in \R: \|\phi_{\theta,1+it}(x)-\phi_{\theta,1+it}(y)\|_1<\alpha\}) \geq 1-\alpha,$$ then
 $\|x-y\|_\theta<\varepsilon$.
 
     \item[(2)] Assuming $X_0$ and $X_1$ are uniformly convex:
 for every $\varepsilon>0$, there exists $\alpha>0$ such that whenever $\|x-y\|_{\theta}<\alpha$ then
 $\gamma_{\theta}^1(T_{x,y}(\varepsilon)) \geq 1-\varepsilon$, where
 $$T_{x,y}(\varepsilon):=\{t \in \R: \|\phi_{\theta,1+it}(x)-\phi_{\theta,1+it}(y)\|_1<\varepsilon\}.$$
\end{itemize}
 \end{lemma} 
 
 \pf

 (1) We know that
 $\|x-y\|_\theta = ||F(\theta)-G(\theta)\|_\theta \leq \|F-G\|_ {\mathcal{F}_\theta^2}$, if $F=F_\theta^x$ and $G=F_\theta^x$. Note also that $F(1+it)=\phi_{\theta,1+it}(x)$ and
 $G(1+it)=\phi_{\theta,1+it}(y)$ for almost all $t \in \R$.
 Using the fact that $F-G$ is bounded by $2$ on $\partial \mathbb{S}$, 
 it then follows from Lemma \ref{border}(2) that there exists $\alpha'>0$ such 
 that  $\|x-y\|_\theta < \varepsilon$ whenever there exists $T$ of measure at least $1 - \frac{\alpha'}{2}$ such that $\|F(1+it)-G(1+it)\|_1<\frac{\alpha'}{2}$ for all $t$ in some $T$. 
 Choose $\alpha=\alpha'/2$.
 
 (2) By Lemma \ref{border}(1), pick $\delta$  small enough so that whenever $\|F\|_{{\mathcal F}_\theta^2} < \delta$, it follows that $\|F(1+it)\|_{1} < \varepsilon$
 except on a set of measure at most $\varepsilon$. Then by Daher's result of uniform continuity of the map $\Gamma_\theta: S_{X_\theta} \rightarrow {\mathcal F}_\theta^2$, we can pick $\alpha>0$ such that
 $$\|F_x^\theta-F_y^\theta\|_{\mathcal{F}_\theta^2} < \delta$$ whenever $\|x-y\|_\theta<\alpha$. 
 \pff
 
 \subsection{Vertical uniformities}
 
 At this point Lemma \ref{55} means that two points $x, y$ of $S_{X_\theta}$
 being close corresponds to the values in $x, y$ of the limit maps ``at level $t$'' being close for ``many'' values of $t$. In order to prove the existence of uniform homeomorphisms between the unit spheres of $X_\theta$ and $X_1$, one can think of some options: a) choosing one appropriate level $t$ and considering the limit map $\phi_{\theta,1+it}$, or b) doing an averaging process of the limit maps over $t$. However in a) there does not seem to be a way of choosing   
 the level $t$ in a uniform way in $x$, and in b) there does not seem to be an averaging process on the sphere $S_{X_1}$ whose result remains in $S_{X_1}$.
 
 We are therefore led to additional assumptions on the scale, which would guarantee that the values of the limit maps
 $\phi_{\theta,1}(x)$ and $\phi_{\theta,1}(y)$ being close (i.e. close at the single level $t=0)$ imply that the values $\phi_{\theta,1+it}(x)$ and $\phi_{\theta,1+it}(y)$ at level $t$ are also close for ``many'' values of $t$. This is achieved most naturally by hypotheses of uniformity of the moduli of the vertical maps on the sphere of $X_s$, for $s$ arbitrarily close to $1$, which motivates the next definitions.

 To simplify notation we shall say that a map $f$ between metric spaces is $(\alpha,\varepsilon)$-continuous if
 $d(f(x),f(y))<\varepsilon$ whenever 
 $d(x,y)<\alpha$. So $f$ is uniformly continuous if for any $\varepsilon>0$ there exists $\alpha>0$ such that $f$
 is $(\alpha,\varepsilon)$-continuous.
 When we say that the $\varepsilon$-modulus of continuity of a family $f_s, s \in S$ is uniform over $S$, we mean that for some $\alpha>0$, $f_s$ is 
 $(\alpha,\varepsilon)$-continuous for all $s \in S$.
 
\begin{defi} \label{def:uniformities} Let $(X_0,X_1)$ be a regular and optimal couple of interpolation spaces. We say that this interpolation scale is
\begin{itemize}
\item[(a)] {\em locally uniform} to the right if there exists a set of $\gamma_{1/2}^1$-positive measure $T \subseteq \R$, such that for each $t \in T$ and each $\varepsilon>0$ the $\varepsilon$-modulus of uniform continuity of vertical maps $\phi_s^t$ is uniform over all $1/2<s<1$.

\item[(b)]  {\em  partially uniform} to the right if  for all $\varepsilon>0$, there exists a set $T_\varepsilon \subseteq \R$ of positive $\gamma_{1/2}^1$-measure, such that the $\varepsilon$-modulus of uniform continuity of vertical maps $\phi_s^t$ is uniform over all $1/2<s<1$ and $t \in T_\varepsilon$.

\item[(c)] {\em  uniform} to the right if  for all $\varepsilon>0$, there exists a set $T_\varepsilon \subseteq \R$ of $\gamma_{1/2}^1$-measure at least $1-\varepsilon$, such that the $\varepsilon$-modulus of uniform continuity of vertical maps $\phi_s^t$ is uniform over all $1/2<s<1$ and $t \in T_\varepsilon$.

\end{itemize}
\end{defi}

It is a small computation that if $X_0$ and $X_1$ are uniformly convex, then we obtain equivalent definitions by replacing $1/2<s<1$
by $\theta_0<s<1$, for any choice of $0<\theta_0<1$. 

It is trivial that (c) $\Rightarrow$ (b), and we also note that  (c) $\Rightarrow$ (a). Indeed if for each $n \in \N$, $T_n$ is associated to $\varepsilon=\varepsilon_n$ in (c), then $T=\cap_{n} T_n$ satisfies (a), and has positive measure provided $\Sigma_{n \in \N}\varepsilon_n$ is small enough. We do not know whether the reverses hold.

\begin{rema} If $X_0$ and $X_1$ are uniformly convex, then items (a), (b), (c) of the previous definition are implied by their versions with the vertical map $\phi_s^t$ substituted by the evaluation map $e_{s,t}$
in $s+it$, i.e. $$e_{s,t}:F \mapsto F(s+it),$$ considered from the sphere of ${\mathcal F}_2^{s}$ into the sphere of $X_s$.
\end{rema}

\pf Indeed $\phi_s^t$ is the composition of $\Gamma_s$ from $S_{X_s}$ to the sphere of 
${\mathcal F}_2^s$ with the above evaluation map
$e_{s,t}$. Since by Proposition \ref{unifunif}, the modulus of continuity of $\Gamma_s$ is uniform over $s$, the result follows.\pff

\begin{prop} Let $(X_0,X_1)$ be a regular compatible couple of uniformly convex spaces.  If the scale admits extrapolation to the right
then it is uniform to the right.
\end{prop}

\pf Because of Fact \ref{compact}, under extrapolation to the right, the maps $\phi_s^t$ have a common modulus of uniform continuity as soon as $0<s<1$ and $t$ belongs to any fixed bounded interval in $\R$.
Indeed

\begin{enumerate}
    \item If $(Y_0, Y_1)$ is a regular pair of uniformly convex spaces then the $\phi_s^t$ are uniformly continuous for $1/2 < s < s_0$ and $t$ in a bounded interval in $\R$, for some $s_0 < 1$, by Fact \ref{compact}.
    \item Under extrapolation, the $\phi_s^t$ of the pair $(X_0, X_1)$ correspond to $\phi_s^t$ of a pair $(X_0, Y)$ with $Y$ uniformly convex.
\end{enumerate}
\pff

\begin{example}
The scale associated to $(\ell_\infty,\ell_1)$ is uniform to the right.
\end{example}

\pf This follows from Example \ref{verticalforlp} for $q=1$.
\pff

\subsection{The limit map as a uniform homeomorphism between spheres}
 
 As previously mentioned, uniformities are candidates for weaker properties than extrapolation, which would imply property (H) or uniform homeomorphism between spheres, through the limit map.
 We first need a technical but important lemma to characterize uniform continuity of the limit map.

 Recall that $\Delta_\theta=\{x \in S_{X_\theta}:
 F_x^{\theta}(s) \overset{w}{\rightarrow} F_x^{\theta}(1) \in S_{X_1}\}$,
 and that $\phi_{\theta,1}$ denotes the map $x \mapsto F_x^{\theta}(1)$ defined on $\Delta_\theta$ with values in $S_{X_1}$.

\begin{lemma}\label{unifor} Assume $(X_0, X_1)$ is an adequate scale of uniformly convex separable spaces.
Assume the moduli of continuity of $\phi_{\theta, s}$ are uniform ``when $s$ tends to $1$'' in the following sense: for any $\varepsilon>0$, there exists $\alpha>0$, such that for all $x,y \in S_{X_\theta}$ with $\|x-y\|_\theta<\alpha$, we have that
$\limsup_{s \rightarrow 1} \|\phi_{\theta,s}(x)-\phi_{\theta,s}(y)\|_{X_0+X_1} <\varepsilon$.

    Then
\begin{enumerate} 
 \item[(1)] for all $x \in S_{X_\theta}$,
    $\phi_{\theta,s}(x)$ converges weakly (necessarily to $\phi_{\theta,1}(x)$).
    \item[(2)] $\phi_{\theta,1}$ is norm to weak continuous from $S_{X_\theta}$ to $X_1$.
\item[(3)]
 if  $\phi_{\theta,1}$ is uniformly continuous on a dense subset of $\Delta_\theta$, then
 $\Delta_\theta=S_{X_\theta}$ and $\phi_{\theta,1}$ is uniformly continuous from $S_{X_\theta}$ to $S_{X_1}$.
 \item[(4)]
 assume furthermore that for any $\varepsilon>0$, there exists $\alpha>0$, such that for all $x,y$ in a dense subset $D_\theta$ of $S_{X_\theta}$, if $\|x-y\|_\theta<\alpha$ then
 $\limsup_{s \rightarrow 1} \|\phi_{\theta,s}(x)-\phi_{\theta,s}(y)\|_{s} <\varepsilon$. 
 Then  $\Delta_\theta=S_{X_\theta}$ and $\phi_{\theta,1}$ is uniformly continuous from $S_{X_\theta}$ to $S_{X_1}$.
 \end{enumerate}

\end{lemma}

\

Note that the hypothesis in (4) holds in particular if  there is $s_0 \in (0, 1)$ such that the modulus of continuity of $\phi_{\theta, s}$ are uniform for $s \in (s_0, 1)$; but this seems to be a stronger hypothesis.

\

\pf (1), (2): 
Recall that for $x \in \Delta_{\theta}$ we have $\phi_{\theta, 1}(x) = w-\lim_{s \rightarrow 1^-} \Gamma_{\theta}(x)(s)$. Take $x \in S_{X_{\theta}}$ and let $x_n$ be a sequence 
of elements of $\Delta_{\theta}$ converging to $x$. By our hypothesis, the modulus is uniform for $s$ close enough to $1$, in the sense that for any $\varepsilon>0$ there is an $\alpha>0$, such that 
 whenever $\|x-x_l\|_\theta<\alpha$, then there exists $s_l$ such that
 $\|\phi_{\theta,s}(x)-\phi_{\theta,s}(x_l)\|_{X_0+X_1} \leq \varepsilon$ for all $s \geq s_l$.
 
Also, since $X_0^*$ and $X_1^*$ are separable, and $(X_0 + X_1)^* = X_0^* \cap X_1^*$ it follows that the weak topology on the unit ball of $X_0 + X_1$ is metrizable. 
To define the metric, let $(\phi_n)$ be a dense sequence in the sphere of $X_0^* \cap X_1^*$ and let $d(x, y) = \sum_n 2^{-n} |\phi_n(x - y)|$.

By the properties of minimal functions we have that $$\|\phi_{\theta,s}(x_n)\|_{X_0 + X_1} \leq \|\phi_{\theta,s}(x_n)\|_{X_s} = 1$$ for every $s$ in the strip.

Fix $\varepsilon>0$. We pick $N_0$ large enough so that whenever $N \geq N_0$, then
 $\|x_N-x\|_\theta \leq \alpha$. There exists $s_N^0$ such that for all $s>s_N^0$, $\|\phi_{\theta,s}(x)-\phi_{\theta,s}(x_N)\|_{X_0+X_1} \leq \varepsilon$.
In particular,
$d(\phi_{\theta,s}(x),\phi_{\theta,s}(x_N)) \leq \varepsilon$ whenever $s>s_N^0$.

By the assumption that $x_N$ is in $\Delta_\theta$, there exists $s_N^1$ (depending on $N$), such that
$d(\phi_{\theta,s}(x_N),\phi_{\theta,1}(x_N)) \leq \varepsilon$ for all $s>s_N^1$.

Thus, we deduce that $$d(\phi_{\theta,1}(x_N)),\phi_{\theta,s}(x)) \leq 2\varepsilon$$ for all $s>s_N:=max(s^0_N,s^1_N)$.

One consequence is that whenever $N,M \geq N_0$, through appropriate $s$ we obtain that
$$d(\phi_{\theta,1}(x_N)),\phi_{\theta,1}(x_M)) \leq 4\varepsilon$$
Hence, $\phi_{\theta,1}(x_n)$ is a $d$-Cauchy sequence and therefore converges weakly to some $y_1$ in $X_1$. Also for any $\varepsilon$ pick $N$ so that $d(\phi_{\theta,1}(x_N),y_1) < \varepsilon$. Then the above gives that 
$d(y_1,\phi_{\theta,s}(x)) \leq 3\varepsilon$ for all $s>s_N$.
Therefore
the $d$-limit $y_1$ of
$\phi_{\theta,1}(x_n)$ is the $d$-limit when $s$ tends to $1$ of
$\phi_{\theta,s}(x)$. This proves (1). This proves also (2): Indeed let $(y_n)_n$ in $S_{X_\theta}$ tending to some $x \in S_{X_\theta}$. Let $(x_n)_n$ be a sequence in $\Delta_\theta$ such that
$\|x_n-y_n\|_\theta$ tends to $0$; by the above we may also assume that $d(\phi_{\theta,1}(x_n),\phi_{\theta,1}(y_n))$ tends to $0$.
Since by (1) $\phi_{\theta,1}(x_n)$ tends weakly to $\phi_{\theta,1}(x)$, it follows that also $\phi_{\theta,1}(y_n)$ tends weakly to $\phi_{\theta,1}(x)$. This proves (2).

(3) In this hypothesis, $\phi_{\theta,1}$ admits a uniformly continuous extension $\tilde{\phi}_{\theta,1}$ to $S_{X_\theta}$ with values in $S_{X_1}$. By the uniqueness of this extension, $\tilde{\phi}_{\theta,1}(x)$ is the limit in norm of $\phi_{\theta,1}(x_n)$ (and also the weak limit) if $x_n \in \Delta_\theta$ tends to $x$. Therefore $\tilde{\phi}_{\theta,1}(x)$ coincides with $\phi_{\theta,1}(x)$, so $\phi_{\theta,1}$ itself is uniformly continuous.

(4) Let $\varepsilon'>0$ be such that $2\beta(\varepsilon')<\varepsilon$ (where $\beta$ is the map considered in Lemma \ref{new}), and let $\alpha'$ associated to $\varepsilon'$ in the hypothesis of assertion (4). Let $x, y$ be arbitrary in $D_\theta$ with $\|x-y\|_{\theta}<\alpha'$.
By Lemma \ref{new} we may pick $x',y'$ in the vertical orbit of $x$ and $y$ respectively, arbitrarily close to $x,y$ and in particular such that they are at distance less than $\alpha'$, and satisfying the condition that $\|\phi_{\theta,1}(x')-\phi_{\theta,1}(y')\|_1 \leq 2\beta(\limsup_s\|\phi_{\theta,s}(x')-\phi_{\theta,s}(y')\|_{s}$). By the hypothesis of present assertion (4) 
we therefore deduce that $\|\phi_{\theta,1}(x')-\phi_{\theta,1}(y')\|_{1} \leq 2\beta(\varepsilon') \leq \varepsilon$.
Since $x',y'$ could be arbitrarily close to $x,y$, let $x^{\prime}_n$ and $y^{\prime}_n$ be sequences of elements of this form converging in norm to $x$ and $y$ respectively. By the norm to weak continuity of $\phi_{\theta,1}$ from (2), we deduce that
$\phi_{\theta,1}(x^{\prime}_n)-\phi_{\theta,1}(y^{\prime}_n)$ converges weakly to $\phi_{\theta,1}(x)-\phi_{\theta,1}(y)$, and therefore that
$\|\phi_{\theta,1}(x)-\phi_{\theta,1}(y)\|_{1} \leq \varepsilon$. So $\phi_{\theta,1}$ is uniformly continuous on $D_\theta$ and (3) applies.
\pff

 \begin{prop}\label{limitlimit}
 Let $(X_0,X_1)$ be a regular compatible couple of uniformly convex spaces. Let $0<\theta<1$. 
 \begin{itemize}
     \item[(1)] Assume the scale is locally uniform to the right. If the scale is finite-dimensional then the limit map $\phi_{\theta,1}$ is injective from $\Delta_\theta$ into $S_{X_1}$. Therefore if the scale has a common $1$-FDD $(Y_n)_{n\in\mathbb N}$, then the limit map $\phi_{\theta,1}$ is injective  from $\Delta_\theta \cap c_{00}(Y_n)$ into $S_{X_1}$.
     \item[(2)]  If the scale is partially uniform to the right, then the limit map $\phi_{\theta,1}$ is  uniformly continuous from $S_{X_\theta}$ into $S_{X_1}$. 
     \item[(3)]  If the scale is uniform to the right, then the limit map $\phi_{\theta,1}$ is  a uniform homemomorphism  from $S_{X_\theta}$ onto a closed subset of  $S_{X_1}$. 
 \end{itemize}

 \end{prop}
 
 \pf (1) To prove injectivity, 
 assume  $\phi_{\theta,1}(x)=\phi_{\theta,1}(y)$. We claim that $\phi_{\theta,1}(\phi_\theta^t(x))=\phi_{\theta,1}(\phi_\theta^t(y))$ for all $t \in T$, where $T$ has positive measure. This means that the functions $F_x^{\theta}$ and $F_y^{\theta}$ coincide almost everywhere on $1+iT$, therefore are equal by Lemma \ref{nonconstant}, and so $x=y$.

To prove the claim,  note by assumption that $(F_x^{\theta}-F_y^{\theta})(s)$ tends to $0$ when $s$ tends to $1$. This convergence is weak in $X_0+X_1$, but under the finite dimensional assumption, all the norms involved are uniformly equivalent, so
$\|(F_x^{\theta}-F_y^{\theta})(s)\|_s$ tends to zero.
By local right uniformity, given any $t \in T$, $\|\phi_s^t(F_x^{\theta}(s))-\phi_s^t(F_y^{\theta}(s))\|_s$ also tends to $0$. But the quantity
$$\phi_s^t(F_x^{\theta}(s))-\phi_s^t(F_y^{\theta}(s))$$
is equal to $$F_x^{\theta}(s+it)-F_y^{\theta}(s+it),$$ which tends to $\phi_{\theta,1}(\phi_\theta^t(x))-\phi_{\theta,1}(\phi_\theta^t(y))$ when $s$ tends to $1$ and when the limits exist, i.e. for $t$ in a set of positive measure. This proves the claim.

\

(2) Fix $\varepsilon >0$. 
By Lemma \ref{new},  for $x,y$ arbitrary in $S_{X_\theta}$, let $U_{x,y} \subset \mathbb{R}$ denote a set with full  $\gamma_{1/2}^1$-measure satisfying:

\medskip
 \emph{Point 1. For any $t \in U_{x,y}$, $$\limsup_{s \rightarrow 1^{-}} \|F_x^\theta(s+it)-
F_y^\theta(s+it)\|_s 
\leq \|F_x^\theta(1+it)-F_y^\theta(1+it)
\|_1.
$$}
 \medskip

Let $T$ be given by the definition of partial uniformity to the right for $\varepsilon$, i.e. $T \subset \mathbb{R}$ is a set of positive $\gamma_{1/2}^1$-measure such that the $\varepsilon$-modulus of uniform continuity of the vertical maps $\phi_s^t$ are uniform over all $1/2 < s < 1$ and $t \in T$.
So let $0<\delta< \gamma_{\frac{1}{2}}^1(-T)$ be such that all vertical maps $\phi_s^t$ are $(\delta,\varepsilon)$-uniform for $1/2<s<1$ and $t \in T$. Summing up:

\medskip
 \emph{Point 2. There is a set $T \subset \mathbb{R}$ of $\gamma_{1/2}^1$-positive measure such that for all $1/2 < s < 1$ and $t \in T$, for any $x,y \in S_{X_s}$, if $\|x - y\|_s < \delta$ then $\|\phi_s^t(x) - \phi_s^t(y)\|_s \leq \varepsilon$.}
 \medskip 

Finally by applying Lemma \ref{55} (2) (using $\gamma_{1/2}^1$ instead of the equivalent measure $\gamma_\theta^1$) we find $\alpha > 0$ satisfying the following point:

\medskip
\emph{Point 3. Whenever $x,y \in S_{X_\theta}$ satisfy $\|x - y\|_{\theta} < \alpha$, there is $S_{x, y} \subset \mathbb{R}$ of $\gamma_{1/2}^1$-measure at least $1 - \delta$ such that for every $t \in S_{x, y}$ we have $\|\phi_{\theta, 1 + it}(x) - \phi_{\theta, 1 + it}(y)\|_1 < \delta$.}
\medskip

So let us now fix $x,y \in S_X$ with $\|x-y\|_\theta<\alpha$.
 Since $U_{x, y}$ has full measure and $\gamma_{\frac{1}{2}}^1(S_{x, y}) \geq 1 - \delta > 1 - \gamma_{\frac{1}{2}}^1(-T)$, there exists $t_{x,y} \in (-T) \cap S_{x, y} \cap U_{x,y}$. If we write $F = F_x^{\theta}$ and $G = F_y^{\theta}$ then we have $\|F(1+it_{x, y}) - G(1 + it_{x,y})\|_1 < \delta$ by Point 3.
 By Point 1, we have that 
 $\|F(s+it_{x,y})-G(s+it_{x,y})\|_s<\delta$ for all $s$ close enough to $1$. Since $F(s)=\phi_{s}^{-t_{x,y}}(F(s+it_{x,y}))$ and similarly for $G$, partial uniformity (Point 2) guarantees that for those $s$,
 $\|F(s)-G(s)\|_s \leq \varepsilon$.
 What we have just proved is slightly stronger than the hypothesis of Lemma \ref{unifor} (4), namely
 \emph{for any $\varepsilon>0$, there exists $\alpha>0$, such that for all $x,y \in S_{X_\theta}$ with $\|x-y\|_\theta<\alpha$,  there exists $s_{xy}<1$ such that $\|\phi_{\theta,s}(x)-\phi_{\theta,s}(y)\|_{s} \leq \varepsilon$ for all $s_{xy} \leq s <1$.}
 
 The conclusion  of Lemma \ref{unifor}(4) then asserts that $\phi_{\theta,1}$ is uniformly continuous from $S_{X_\theta}$ to $S_{X_1}$. 

\

 (3)  We need to prove that the inverse map is well defined and uniformly continuous. 
So fix $\varepsilon>0$, and let $\alpha>0$ be such that $\alpha':=\max(\alpha,2\beta(\alpha))$ is associated to $\varepsilon$ in Lemma \ref{55}(1), in which we replaced $\gamma_\theta^1$ by $\gamma_{1/2}^1$; explicitely: whenever $\gamma_{1/2}^1(T') \geq 1-\alpha'$ and $\|\phi_{\theta,1+it}(x)-\phi_{\theta,1+it}(y)\|_1<\alpha'$ for all $t \in T'$, it follows that $\|x-y\|_\theta<\varepsilon$. Using right uniformity, let $\delta>0$ and $T$ of $\gamma_{1/2}^1$-measure at least $1-\alpha/2$ be such that the maps $\phi_s^t$ are $(\delta,\alpha)$-continuous for all $1/2<s<1$ and $t \in T$.

 Let $x,y$ in $S_{X_\theta}$ be such that $\|\phi_{\theta,1}(x)-\phi_{\theta,1}(y)\|_1 <\delta$. By continuity of $\phi_{\theta,1}$ from (2), for $t>0$ in some small enough interval $I$ around $0$, the same inequality holds for $x'=\phi_\theta^t(x)$ and 
$y'=\phi_\theta^t(y)$,  i.e.
$$\|\phi_{\theta,1+it}(x)-\phi_{\theta,1+it}(y)\|_1 <\delta.$$
By Lemma \ref{new} we may replace $I$ by a full measure subset of it  (still called $I$) so that for all $t \in I$,
$$\|\phi_{\theta,s+it}(x)-\phi_{\theta,s+it}(y)\|_s <\delta,$$ for $s$ close enough to $1$. Pick $t_0 \in I$ small enough so that
$t_0+T$ has $\gamma_{1/2}^1$-measure at least $1-\alpha \geq 1-\alpha'$, and $s_0<1$ such that
$$\|\phi_{\theta,s+it_0}(x)-\phi_{\theta,s+it_0}(y)\|_1 <\delta,$$
for all $s \geq s_0$.
By uniformity we have that
$$\|\phi_{\theta,s+i(t_0+t)}(x)-\phi_{\theta,s+i(t_0+t)}(y)\|_s <\alpha$$ for all $t \in T$ and $s \geq s_0$. I.e.
$$\|\phi_{\theta,s+it}(x)-\phi_{\theta,s+it}(y)\|_s <\alpha$$
for all $t \in t_0+T$. Taking the limit when $s$ tends to $1$ and applying Lemma \ref{new}, we have
$$\|\phi_{\theta,1+it}(x)-\phi_{\theta,1+it}(y)\|_1 <2\beta(\alpha) \leq \alpha'$$
a.e. on $t_0+T$.
Since $t_0+T$ has $\gamma_{1/2}^1$-measure at least $1-\alpha'$, Lemma \ref{55} (1) and the choice of $\alpha'$ imply  that $\|x-y\|_\theta$ is smaller than $\varepsilon$.

Summing up, we have proved that $\|x-y\|_\theta<\varepsilon$ for $x, y \in S_{X_\theta}$, as long as $\|\phi_{\theta,1}(x)-\phi_{\theta,1}(y)\|_1<\delta$.
Therefore $\phi_{\theta,1}$ is injective and its inverse is uniformly continous.
\pff

\begin{prop}\label{limitlimitgeneral}
 Let $(X_0,X_1)$ be a regular compatible couple of uniformly convex spaces which is partially uniform to the right. Let $0<\theta<1$. 
 \begin{itemize}
     \item[(1)] If $X_0$ and $X_1$ are finite dimensional, 
 then the limit map $\phi_{\theta,1}$ is a homeomorphism between $S_{X_\theta}$ and $S_{X_1}$.
     \item[(2)] If the scale has common $1$-FDD $(Y_n)_{n\in \mathbb N}$, then $\phi_{\theta,1}$ is a uniformly continuous map from $S_{X_\theta}$ into $S_{X_1}$, whose restrictions to the initial segments $[Y_1,\ldots,Y_n]$ are homeomorphisms.
     \item[(3)] If the scale has common $1$-FDD $(Y_n)_{n\in\mathbb N}$, and is  uniform to the right, then $\phi_{\theta,1}$ is a uniform homeomorphism between $S_{X_\theta}$ and $S_{X_1}$ leaving each $[Y_1,\ldots,Y_n]$ invariant.
 \end{itemize}

 \end{prop}
 
 \pf
(1) By Proposition \ref{limitlimit} (2), we know that the limit map defines a continuous map of $S_{X_\theta}$ onto a closed subset of $S_{X_1}$. Furthermore we know by Proposition \ref{limitlimit} (1) that it is injective on $S_{X_\theta}$ (because it is injective on $\Delta_\theta$, and because the proof of Proposition \ref{limitlimit}(2) indicates that the hypothesis of Lemma \ref{unifor}(4) holds, and therefore that $\Delta_\theta=S_{X_\theta}$). Therefore the limit map defines a homeomorphism of $S_{X_\theta}$ onto its image. Since the scale is regular, $X_0$ and $X_1$ have the same dimension.   
Therefore, the image of $S_{X_{\theta}}$ cannot be a proper subset of $S_{X_1}$, since a finite dimensional sphere is not homeomorphic to any proper subspace (this is a direct consequence of Borsuk-Ulam theorem).

 \
 
 (2) Note that (1) and Proposition \ref{limitlimit} (2) give that the limit map is uniformly continuous,  and that its restrictions to the initial segments $[Y_1,\ldots,Y_n]$ are homeomorphisms.

\

(3) follows from (2) and the fact the inverse map exists and is uniformly continuous, Proposition \ref{limitlimit}(3). Since the image of $S_{X_\theta}$ is closed and since it contains $S_{X_1} \cap c_{00}(Y_n)$, it must be equal to $S_{X_1}$, concluding the argument. 
\pff

\

We obtain our  main result as a corollary: under uniformity of the scale we achieve uniform homemomorphism of the spheres of the spaces in the interior and in the boundary of the domain.

\begin{coro}  \label{cor:main_result} 
Let $X$ be a uniformly convex complex space with a $1$-FDD. 
Consider the interpolation scale $(\overline{X}^* ,X)$. If this scale is  uniform to the right, then $S_X$ is uniformly homeomorphic to $S_{\ell_2}$.
\end{coro}

\begin{conj}
It seems natural to wonder whether Proposition \ref{limitlimitgeneral} can be modified to get uniform continuous maps from $S_{X_1}$ to $S_{X_\theta}$ so that initial segments of the FDD are homeomorphisms. This should yield property (H) for $X$ with an FDD and the property that $(\overline{X}^*,X)_{1/2}=\ell_2$.
\end{conj}

\newpage

\newpage

\end{document}